\documentclass[11pt, a4paper]{amsart}

\usepackage{amsfonts,amsmath,amssymb, amscd}

\usepackage[all]{xy}

\newtheorem{theorem}{Theorem}[section]
\newtheorem{lemma}[theorem]{Lemma}
\newtheorem{prop}[theorem]{Proposition}
\newtheorem{corollary}[theorem]{Corollary}

\theoremstyle{definition}
\newtheorem{definition}[theorem]{Definition}
\newtheorem{rem}[theorem]{Remark}

\newtheorem{example}[theorem]{Example}
\newtheorem{problem}[theorem]{Problem}

\usepackage{color}

\newcommand\pf{\begin{proof}}
\newcommand\epf{\end{proof}}

\newcommand{\ot}{\otimes}
\newcommand{\mr}{\mathrm}
\newcommand{\ms}{\mathsf}
\newcommand{\op}{\operatorname}

\renewcommand{\d}{\delta}
\newcommand{\dK}{\delta\text{-}K}

\renewcommand{\c}{\mathbb{C}}
\newcommand{\C}{\mathbb{C}(t)}
\newcommand{\G}{\mathsf{G}}
\newcommand{\g}{\mathfrak{g}}
\newcommand{\X}{\mathsf{X}}

\newcommand{\Lie}{\operatorname{Lie}}
\renewcommand{\l}{\ell}

\newcommand{\rcosmash}{\mathop{\raisebox{0.2ex}{\makebox[0.92em][l]{${\scriptstyle\blacktriangleright\mathrel{\mkern-4mu}<}$}}}}

\numberwithin{equation}{section}

\title[Twisted forms of differential Lie algebras]
{Twisted forms of differential Lie algebras over $\C$ associated with\\ complex simple Lie algebras}

\author[A.~Masuoka]{Akira~Masuoka}
\address{Akira Masuoka,
Institute of Mathematics, 
University of Tsukuba, 
Ibaraki 305-8571, Japan}
\email{akira@math.tsukuba.ac.jp}

\author[Y.~Shimada]{Yuta~Shimada}
\address{Yuta Shimada,
Graduate School of Pure and Applied Sciences, 
University of Tsukuba, 
Ibaraki 305-8571, Japan}
\email{shimada@math.tsukuba.ac.jp}

\dedicatory{Dedicated to Professor Jun Morita on the occasion of his retirement}

\begin{document}

\begin{abstract}
Discussed here is descent theory in the differential context where everything is equipped with a differential operator. 
To answer a question personally posed by A.~Pianzola, we determine all twisted forms of the differential 
Lie algebras over $\C$ associated with complex simple Lie algebras. Hopf-Galois Theory, a ring-theoretic counterpart
of theory of torsors for group schemes, plays a role when we grasp the above-mentioned twisted forms from torsors. 
\end{abstract}

\maketitle


\medskip


\section{Introduction---problem and answer}\label{sec:intro}

Rings and algebras are supposed to be associative and containing $1$, and their morphisms are supposed to send $1$ to $1$. 
Moreover, rings, algebras and Hopf algebras are assumed to be commutative, unless otherwise stated.

We let $\d$ mean ``differential" and use the symbol $\d$ to indicate differential operators in general. 
A $\d$-\emph{ring} is thus a (commutative) ring $R$ equipped with an additive operator $\d : R \to R$
satisfying the Leibniz rule 
\[
\d(xy)=(\d x)y+ x(\d y)
\] 
for all $x, y \in R$.
It is called a $\d$-\emph{field}
if the ring is a field. The rational function field $\mathbb{C}(t)$ in one variable is regarded as a $\d$-field with respect to the 
standard operator such that $\d t=1$ and $\d c=0$ for every $c \in \mathbb{C}$.
The field
\begin{equation*}\label{eq:constants}
C_{\C}= \{ \, a \in \C \mid \d a=0 \, \}
\end{equation*}
of constants is $\mathbb{C}$. 
A $\d$-$\mathbb{C}(t)$-\emph{Lie algebra}
is a Lie algebra $\mathfrak{g}$ over $\mathbb{C}(t)$ which is equipped with an additive operator $\d : \mathfrak{g}\to \mathfrak{g}$
such that
\[ \d(a X)= (\d a)X+ a(\d X),\quad \d[X,Y]=[\d X, Y]+ [X,\d Y] \]
for all $a \in \C$ and $X,Y \in \mathfrak{g}$. In the same way a $\d$-$R$-\emph{Lie algebra} is defined for any $\d$-ring $R$.
Let $\mathfrak{g}$ be a $\d$-$\mathbb{C}(t)$-Lie algebra. Given a $\d$-ring map $\mathbb{C}(t)\to R$ (that is, a ring map
preserving the $\d$-operator), the base extension $\mathfrak{g} \ot_{\mathbb{C}(t)} R$ of $\mathfrak{g}$ is naturally a 
$\d$-$R$-Lie algebra. A \emph{twisted form} of $\mathfrak{g}$ is a $\d$-$\mathbb{C}(t)$-Lie algebra $\mathfrak{f}$ such that
\[\mathfrak{g} \ot_{\mathbb{C}(t)} R \simeq \mathfrak{f} \ot_{\mathbb{C}(t)} R \quad \text{as}\ \d\text{-}R\text{-Lie algebras} \]
for some $\d$-ring map $\mathbb{C}(t)\to R$ with $R\ne 0$.

Let $n \ge 2$. We can and do regard the $\mathbb{C}(t)$-Lie algebra $\mathfrak{sl}_n(\mathbb{C}(t))$ which consists of all traceless matrices 
$X=\big(x_{ij}\big)$ with $x_{ij}\in \mathbb{C}(t)$, as a $\d$-$\mathbb{C}(t)$-Lie algebra
with respect to the entry-wise $\d$-operator $\d \big(x_{ij}\big) := \big(\d x_{ij}\big)$.\
A.~Pianzola posed personally to the first-named author, a problem which is 
specialized to the following. 

\begin{problem}
Describe all twisted forms of the $\d$-$\mathbb{C}(t)$-Lie algebra $\mathfrak{sl}_n(\mathbb{C}(t))$. 
\end{problem}

Clearly, to generalize the problem, one can replace $\mathfrak{sl}_n(\c)$ with a complex simple Lie algebra 
$\g_0$ (of finite dimension), and regard $\C$-Lie algebra
\begin{equation}\label{eq:Lie_functor}
\g_0(\C)=\g_0 \otimes_\c \C
\end{equation}
as a $\d$-$\C$-Lie algebra with the $\d$ operating on the tensor factor
$\C$.
The notation \eqref{eq:Lie_functor} is used since $\g_0$ is seen to give the functor $R\mapsto \g_0 \otimes_\c R$,
and $\g_0 \otimes_\c \C$ is then its value. 
Pianzola's problem is precisely the following one, 
which we are going to solve. 

\begin{problem}
Given a complex simple Lie algebra $\g_0$, describe all twisted forms of $\g_0(\C)$. 
\end{problem}

Our interest in this problem or in studying twisted forms of Lie algebras 
in the differential context comes from differential Galois Theory, in which examples of such twisted forms naturally arise. 
In fact, given a homogeneous linear differential equation
with coefficients in a differential field, say $\C$, 
we have the \emph{Galois group} of the equation, $\mathsf{G}_0$, which is an affine algebraic $\c$-group. 
In addition we naturally have an affine algebraic
``differential" $\C$-group $\mathsf{G}$, called the
\emph{intrinsic Galois group}
(or \emph{Katz group}); see \cite{Katz}. 
One sees that the Lie algebra $\Lie(\mathsf{G})$ of $\mathsf{G}$ is 
a $\d$-$\C$-Lie algebra
which is a twisted form of the $\d$-$\C$-Lie algebra $\Lie(\mathsf{G}_0)\otimes_{\c}\C$,
where $\Lie(\mathsf{G}_0)$ is the (complex) Lie algebra of $\mathsf{G}_0$.

Returning to the situation before,
recall that the complex simple Lie algebras are classified, 
labeled by their root systems
\[
\mr{A}_\l~(\l \ge 1),\ \mr{B}_\l~(\l \ge 2),\ \mr{C}_\l~(\l \ge 3),\ \mr{D}_\l~(\l \ge 4),\ \mr{E}_6,\ \mr{E}_7,\ \mr{E}_8,\ \mr{F}_4,\ \mr{G}_2. 
\] 
See \cite[Chapter IV]{J}, for example.
Let $\g_0$ be a complex simple Lie algebra, and
let $\Gamma$ denote the automorphism group of the associated Dynkin diagram. 
Explicitly, the group is
\begin{equation}\label{eq:Gamma}
\Gamma =
\begin{cases}
\{ 1 \} & \text{type}\ \mr{A}_1,\ \mr{B}_\l~(\l \ge 2),\ \mr{C}_\l~(\l \ge 3),\ \mr{E}_7,\ \mr{E}_8,\ \mr{F}_4\ \text{or}\ \mr{G}_2;\\
\mathbb{Z}_2 & \text{type}\ \mr{A}_{\ell}\ (\ell \ge 2),\ \mr{D}_{\ell}\ (\ell \ge 5)\ \, \text{or}\ \, \mr{E}_6; \\
\mathfrak{S}_3 & \text{type}\ \mr{D}_4
\end{cases}
\end{equation}
according to the type of $\g_0$; see \cite[Table 3 on Page 298]{OV}. 
Here and in what follows $\mathbb{Z}_n$ denotes the cyclic group of order $n$. 
In addition, $\mathfrak{S}_3$ denotes the symmetric group
of degree 3.
The action by $\Gamma$ naturally (up to conjugation) gives rise to automorphisms of $\g_0$, which forms a group naturally identified with
the group $\ms{Out}(\g_0)$ of outer-automorphisms of $\g_0$. 

Roughly speaking, our answer, Theorem \ref{thm1}, to the problem tells
that all non-trivial twisted forms are obtained by the Galois descent (see \cite[Section 18]{KMRT})
for which $\Gamma$ (and its subgroups for type $\mathrm{D}_4$)
act as Galois groups. 
To make a precise statement we introduce below the notion of being \emph{quasi-isomorphic}.  

\begin{lemma}\label{lem:dadD}
If $\mathfrak{g}=(\mathfrak{g}, \d)$ is a $\d$-$R$-Lie algebra, then for any element $D \in \mathfrak{g}$,
\[ \d + \mathsf{ad}(D) : \mathfrak{g}\to \mathfrak{g},\ X \mapsto \d X +[D,X] \]
is a $\d$-operator with which $\mathfrak{g}$ is again a $\d$-$R$-Lie algebra. 
\end{lemma}

Indeed, one sees, more generally, that for any $R$-linear derivation $\mathfrak{D} : \mathfrak{g}\to \mathfrak{g}$, 
$(\mathfrak{g}, \d+\mathfrak{D})$ is a $\d$-$R$-Lie algebra. 
Note that the inner derivation $\mathsf{ad}(D)$ above is $R$-linear. 

\begin{definition}\label{def:quasi-isom}
Let $R$ be as above, 
We say that two $\d$-$R$-Lie algebras 
$\mathfrak{g}=(\mathfrak{g}, \d)$ and $\mathfrak{g}'=(\mathfrak{g}', \d')$
are \emph{quasi-isomorphic},
if there is an element $D \in \mathfrak{g}$ such that
\[ (\mathfrak{g}, \d+\mathsf{ad}(D))\simeq (\mathfrak{g}', \d') \quad \text{as}\ \d\text{-}R\text{-Lie algebras}. \]
\end{definition}

The condition is equivalent to saying that there is an element $D' \in \mathfrak{g}'$ such that
$(\mathfrak{g}, \d)\simeq (\mathfrak{g}', \d'+\mathsf{ad}(D'))$, as is easily seen. 
It follows that the quasi-isomorphism gives an equivalence
relation among all $\d$-$R$-Lie algebras.

\begin{theorem}\label{thm1} 
Suppose that $\g_0$ is a complex simple Lie algebra, and let $\Gamma$ be the automorphism group of the associated
Dynkin diagram. Then a
$\d$-$\C$-Lie algebra is a twisted form  of 
$\g_0(\C)$ if and only if it is quasi-isomorphic to one of those listed below, according to the case 
$\Gamma=\{1\}$, $\mathbb{Z}_2$ or $\mathfrak{S}_3$; see \eqref{eq:Gamma}. 
\begin{itemize}
\item[\textup{(1)}] Case $\Gamma=\{1\}$\textup{:}\ $\g_0(\C)$\textup{;}

\medskip

\item[\textup{(2)}] Case $\Gamma=\mathbb{Z}_2$\textup{:}~~\textup{(i)}~$\g_0(\C)$\textup{;}
\begin{itemize}
\item[\textup{(ii)}] $\g_0(L)^{\Gamma}$, where $L/\C$ is
a quadratic field extension.
\end{itemize}

\medskip

\item[\textup{(3)}] Case $\Gamma=\mathfrak{S}_3$\textup{:}~~\textup{(i)}~$\g_0(\C)$\textup{;}
\begin{itemize}
\item[\textup{(ii)}] $\g_0(L)^{\mathbb{Z}_2}$, where $L/\C$ is a quadratic field extension\textup{;}
\item[\textup{(iii)}] $\g_0(L)^{\mathbb{Z}_3}$, where $L/\C$ is a cubic Galois extension\textup{;}
\item[\textup{(iv)}] $\g_0(L)^{\Gamma}$, where $L/\C$ is
a Galois extension of fields with Galois group $\Gamma(=\mathfrak{S}_3)$. 
\end{itemize}
\end{itemize}
\end{theorem}

We should immediately add some explanations about the statement above. 
First, any finite field extension $L/\C$ uniquely turns into an extension of $\d$-fields, 
whence $\g_0(L)$ turns into a $\d$-$L$-Lie algebra. 
Second, 
in (ii) of (2) and (iv) of (3) above, the group $\Gamma$ is supposed to act diagonally on $\g_0(L)=\g_0\otimes_{\c} L$,
as outer-automorphisms on $\g_0$, and as the Galois group on $L$. 
In addition, $\g_0(L)^{\Gamma}$ denotes the $\Gamma$-invariants
in $\g_0(L)$, which is  in fact a $\d$-$\C$-Lie algebra by Galois descent; see Section \ref{subsec:Galois_descent}. 
Third, in (ii) of (3), we choose arbitrarily an order 2 subgroup
$\mathbb{Z}_2$ of $\Gamma\, (=\mathfrak{S}_3)$, and let it act on $\g_0$ by restriction.
The $\d$-$\C$-Lie algebra
$\g_0(L)^{\mathbb{Z}_2}$ which results in the same way as above
does not depend (up to isomorphism) on the choice since the order $2$ subgroups are conjugate to each other;
on the other hand it may depend on $L$. 
Finally, in (iii) of (3), 
we suppose that $\mathbb{Z}_3$ is the unique order 3 subgroup of $\Gamma\, (=\mathfrak{S}_3)$,
and let it act on $\g_0$ by restriction, again. 
We add the following remark: there exist infinitely many quadratic and cubic Galois extensions over $\C$, as is easily seen, 
while the existence of a Galois extension over $\C$ with Galois group $\mathfrak{S}_3$ will be ensured 
by Example \ref{example:Galois_extension}.

The theorem will be proved in the final Section \ref{sec4}, which contains as well,
explicit descriptions (see Section \ref{subsec:explicit_form})
of the non-trivial $\d$-$\C$-Lie algebras listed in (ii) of (2) and (ii)-(iv) of (3). 
The preceding two sections provide preliminaries, some of which are beyond what 
will be needed, but are of interest by themselves. 
Section \ref{sec2} presents
descent theory in the differential context; Section \ref{sec3} prepares technical tools mainly from
Hopf-Galois Theory, which is a ring-theoretic counterpart of theory of torsors for group schemes. 
In particular, Schauenburg's bi-Galois Theory \cite{Sc} will play a role 
in two stages (see Sections \ref{subsec:interpret} and \ref{subsec:Galois_descent}),
when we grasp the twisted forms in question from $\d$-torsors.


\section{Descent theory in differential context}\label{sec2}


\subsection{$\d$-$R$-Objects}
Let $R$ be a $\d$-ring. 
A $\d$-$R$-\emph{module} is an $R$-module $M$ equipped with an additive operator 
$\d : M \to M$
satisfying 
\[
\d(xm)=(\d x)m+ x(\d m),\quad x\in R,\ m \in M. 
\]
All $\d$-$R$-modules form a symmetric tensor category $\d$-$R$-$\mathtt{Modules}$ with respect to the 
tensor product $M_1 \ot_R M_2$, 
the unit object $R$ and the obvious symmetry 
\[
M_1 \ot_R M_2\to M_2\ot_R M_1,\quad m_1\ot m_2 \mapsto m_2 \ot m_1. 
\]
The $\d$-operator on $M_1\ot_R M_2$ is given by 
\[
\d(m_1\ot m_2)=\d m_1\ot m_2+m_1\ot \d m_2. 
\]

The notion of $\d$-$R$-Lie algebra defined in the previous section is precisely a Lie algebra in the category $\d$-$R$-$\mathtt{Modules}$.
In general, any linear object, such as algebra or Hopf algebra, in $\d$-$R\text{-}\mathtt{Modules}$ is called
a $\d$-$R$-\emph{object}, so as $\d$-$R$-\emph{algebra} or $\d$-$R$-\emph{Hopf algebra};
important is the fact that the structure is defined by morphisms of $\d$-$R$-$\mathtt{Modules}$ between tensor powers of the object. 
Given a $\d$-$R$-algebra $S$, we have the base-extension functor
\begin{equation*}\label{eq:base_ext_functor}
\ot_R S : \d\text{-}R\text{-}\mathtt{Modules}\to \d\text{-}S\text{-}\mathtt{Modules}, 
\end{equation*}
which induces base-extension functors for linear objects such as above. 

We are concerned with descent theory (see \cite{W}, for example) in differential context. To make this clearer, 
let us fix a $\d$-$R$-object $A$. A $\d$-$R$-object $B$ is called an $S/R$-\emph{form} of $A$, 
or a \emph{twisted form} of $A$
\emph{split by} $S$, if 
$S$ is a $\d$-$R$-algebra such that
\begin{itemize}
\item[(i)] $S$ is faithfully flat as an $R$-algebra, and
\item[(ii)] $A \ot_R S \simeq B \ot_R S$ as $\d$-$S$-objects.
\end{itemize}
A $\d$-$R$-object $B$ is called a \emph{twisted form} of $A$, if there exists a $\d$-$R$-algebra $S$ which satisfies (i) and (ii) above.

The $\d$-\emph{automorphism-group functor} of $A$ is the functor
\begin{equation}\label{eq:Aut}
\mathsf{Aut}_{\d}(A) : \d\text{-}R\text{-}\mathtt{Algebras}\to \mathtt{Groups},\quad T \mapsto \op{Aut}_{\d\text{-}T}(A\ot_R T) 
\end{equation}
from the category $\d\text{-}R\text{-}\mathtt{Algebras}$ of $\d$-$R$-algebras to the category 
$\mathtt{Groups}$ of groups, which associates to each  $\d$-$R$-algebra $T$, the automorphism-group
${\op{Aut}_{\d\text{-}T}(A \ot_R T)}$ of the $\d$-$T$-object $A \ot_R T$. 
When constructing the 1st Amitsur cohomology (pointed) set as in \cite[Section 17.6]{W}, replace faithfully flat
homomorphisms of rings and automorphism-group functors with our $R \to S$ (satisfying (i) above)
and $\mathsf{Aut}_{\d}(A)$, respectively. The resulting
differential analogue is denoted by
\begin{equation}\label{eq:H1}
H^1_{\d}(S/R, \mathsf{Aut}_{\d}(A)).
\end{equation}
This is seen to classify the $S/R$-forms of $A$; to be more precise, 
there is a natural bijection from 
$H^1_{\d}(S/R, \mathsf{Aut}_{\d}(A))$ to the set of 
all $\d\text{-}R$-isomorphism classes of the $S/R$-forms.
An important consequence is: if $A'$ is another $\d$-$R$-object 
of some distinct kind, which has the $\d$-automorphism-group functor isomorphic to $\mathsf{Aut}_{\d}(A)$, 
then there is a one-to-one correspondence (up to isomorphism) 
between the $S/R$-forms of $A$ and 
$S/R$-forms of $A'$. 

We remark that for any functor $\G : \d\text{-}R\text{-}\mathtt{Algebras}\to \mathtt{Groups}$, the cohomology
set $H^1_{\d}(S/R, \G)$ is defined just as the one in \eqref{eq:H1}. The set will appear in what follows
(see \eqref{eq:H1KG}) only when $\G$ is
representable, and turns out, indeed, to be an automorphism-group functor. 

\begin{rem}\label{rem:right_left}
We have used so far the \emph{base-on-right} notation $A \ot_R T$ which denotes the extended base on the right;
it seemingly looks nicer than the \emph{base-on-left} notation $T\ot_R A$. But we may and do (when it is natural)
use the latter notation. 
\end{rem}

\begin{rem}
In this Remark the differential structure is ignored.
\begin{itemize}
\item[(1)]
We would like to clarify our use of the term (twisted) form. 
Consider a faithfully flat homomorphism $R\to S$ of rings and an $R$-object $A$.
If $R$ and $S$ are fields, it is more conventional (as one of the referees
mentioned) for experts in representation theory or physicists to speak of $R$-form
of the $S$-object $A \otimes S$. For example,
$\mathfrak{so}_3(\mathbb{R})$
is a real form of $\mathfrak{sl}_2(\c)$. Our terminology is that 
$\mathfrak{so}_3(\mathbb{R})$ is a
twisted form, or a $\c/\mathbb{R}$-form of $\mathfrak{sl}_2(\mathbb{R})$.  
The terminology we have
chosen, namely an $S/R$-form or twisted form of the $R$-object $A$, is familiar
in number theory and  algebraic geometry, and it is also the standard
terminology in Grothendieck's  descent theory. The base $R$ is fixed, and
the $S$ can vary.
\item[(2)]
The referee suggested the authors to add the articles 
\cite{B1}--\cite{BR2}, \cite{Ser1}
and \cite{Ser2}
to the References.
Those articles classify real forms of \emph{loop algebras}, which refer to the complex
Lie superalgebras of the form
$\g_0(\c[t,t^{-1}])$, where $\g_0$ is a complex simple 
Lie superalgebra of finite dimension. Moreover, they associate with those forms 
affine Kac-Moody superalgebras. 
Notable is Serganova's method (see \cite{Ser2})
based on the classification of real forms of the complexified Lie algebra
of vector fields on a circle.
\item[(3)]
Pianzola introduced to the authors his \cite{P};
it characterizes 
affine Kac-Moody Lie algebras as twisted forms (in the \'{e}tale sense)
of
the untwisted algebras. This paper likely contains ideas
which may be used to replace some of the argument of Section \ref{subsec:Galois_descent}
below by appealing to Galois cohomology.  
\end{itemize}
\end{rem}


\section{Affine $\dK$-groups, their Lie algebras and torsors}\label{sec3}

In this section $K$ denotes a $\d$-field. 
We assume that the characteristic $\op{char} K$ of $K$ is zero. 


\subsection{Affine $\dK$-groups and their Lie algebras}
An \emph{affine $\dK$-group scheme} is by definition 
a representable functor 
\[
\G : \dK\text{-}\mathtt{Algebras} \to \mathtt{Groups}\quad \textup{(see \eqref{eq:Aut})}; 
\]
this will be called an \emph{affine 
$\dK$-group} for short. 
Such a functor $\G$ is uniquely represented by a $\dK$-Hopf algebra, say $H$, and is presented
so as $\G=\op{Spec}_{\d}(H)$ or $\op{Spec}_{\dK}(H)$.  We say that $\G$ is \emph{algebraic}, or it is an \emph{affine algebraic $\dK$-group},
 if $H$ is finitely generated as a $K$-algebra. 
If one forgets $\d$, then $\G=\op{Spec}(H)$ is an affine $K$-group, which has the Lie algebra 
\begin{equation}\label{eq:LieG}
\op{Lie}(\G)=\op{Der}_{\epsilon}(H, K). 
\end{equation}
Recall that this consists of all $K$-linear maps $D : H \to K$ that satisfy 
\[
D(ab)=D(a)\epsilon(b)+\epsilon(a)D(b),\quad a, b \in H, 
\]
where $\epsilon : H \to K$ is the counit of $H$. This is in fact a $\dK$-Lie algebra with respect to the operator defined by
\[ (\d D)(a):=\d(Da)-D(\d a),\quad D \in \op{Lie}(\G),\ a \in H . \]
Note that the canonical pairing $H\ot_K \op{Lie}(\G) \to K$ is a morphism in
$\dK$-$\mathtt{Modules}$. 
We have $\dim_K(\op{Lie}(\G))<\infty$, if $\G$ is algebraic. 

\begin{rem}\label{rem:algebraic}
The notion of being ``algebraic" defined above would be rather restricted for those who would like to work intensively in differential algebra. 
It should be distinguished from the more natural (for those above) notion of being ``$\d$-algebraic", 
which will be discussed briefly in Section \ref{subsec:d-algebraic}, being less crucial for our purpose though. 
\end{rem}


\subsection{$\dK$-Torsors and Galois $\dK$-algebras}
An \emph{affine $\dK$-scheme} is by definition a representable set-valued functor
\[
\X : \dK\text{-}\mathtt{Algebras}\to \mathtt{Sets}.
\]
It is uniquely represented by a $\dK$-algebra, say $A$, being
presented so as $\X= \op{Spec}_{\d}(A)$; it is said to be \emph{algebraic} if $A$ is finitely generated as a $K$-algebra.
The category of affine $\dK$-schemes, whose morphisms are natural transformations, has direct products. 
The direct product $\X_1 \times \X_2$ of two affine $\dK$-schemes $\X_i= \op{Spec}_{\d}(A_i)$, $i=1,2$,
is represented by $A_1 \ot_K A_2$. The notion of \emph{group object} of the category is
naturally defined, and such an object is precisely an affine $\dK$-group. 
Given an affine $\dK$-group $\G=\op{Spec}_{\d}(H)$, the notion of right (or left) $\G$-equivariant objects 
is defined, as well. Such an object is called a \emph{right} (or \emph{left}) 
\emph{$\G$-equivariant $\dK$-scheme}. Giving such a $\dK$-scheme
$\X= \op{Spec}_{\d}(A)$ is the same as giving a \emph{right} (or \emph{left}) \emph{$H$-comodule $\dK$-algebra};
it is an object $A$ in $\dK\text{-}\mathtt{Algebras}$ equipped with a morphism $A \to A\ot_K H$ (or
$A \to H\ot_K A$) in the category which satisfy the co-associativity and the counit property. 
Obviously, $\G$ itself is $\G$-equivariant on both sides. 

Let $R$ be a $\dK$-algebra. An affine $\dK$-group or (equivariant or ordinary) $\dK$-scheme 
$\X= \op{Spec}_{\d}(A)$ has the base change $\X_R=\op{Spec}_{\d\text{-}R}(A\ot_K R)$; it is by definition 
the functor $T \mapsto \X(T)$ defined on $\d\text{-}R\text{-}\mathtt{Algebras}$, where
each $T \in \d\text{-}R\text{-}\mathtt{Algebras}$ is regarded naturally as a $\dK$-algebra.
We can discuss twisted forms of $\X$; it is the same as discussing twisted forms of $A$. 

Let $\G=\op{Spec}_{\d}(H)$ be an affine $\dK$-group. A twisted form of 
the right $\G$-equivariant $\dK$-scheme $\G$ is called a \emph{right} $\dK$-\emph{torsor} for $\G$.
To be explicit it is a right $\G$-equivariant $\dK$-scheme $\X$ such that $\X_R \simeq \G_R$
as right $\G$-equivariant $\d\text{-}R$-schemes for some non-zero $\dK$-algebra $R$. 
Such an $\X$ is uniquely represented by a right $H$-comodule $\dK$-algebra $B$ which 
is a twisted form of $H$. Such a twisted form $B$ is characterized as a \emph{right $H$-Galois $\dK$-algebra}
\cite[Section 8.1]{Mo}; it is by definition a non-zero right $H$-comodule $\dK$-algebra $B$ such that 
the $\dK$-algebra map
\begin{equation}\label{eq:rho}
\tilde{\rho} : B \ot_K B \to B \ot_K H,\quad \tilde{\rho}(b \otimes c)=b \rho(c)
\end{equation}
is an isomorphism. Here and in what follows, $\rho : B \to B \ot_K H$ denotes the structure map. 
Note that $\tilde{\rho}$ is a $\d\text{-}B$-algebra isomorphism 
(with the base-on-left notation, see Remark \ref{rem:right_left}), 
and $B$ is split by $B$ itself. 

The analogous notions of \emph{left $\dK$-torsors for $\G$} and of \emph{left $H$-Galois $\dK$-algebras}
are defined in the obvious manner, and those two are in one-to-one correspondence. 


\subsection{Affine $\d$-algebraic $\dK$-groups}\label{subsec:d-algebraic}
An $\dK$-algebra $A$ is said to be \emph{$\d$-finitely generated} if it is generated as a $K$-algebra by finitely
many elements $a_1,\dots,a_n$ together with their iterated differentials $\d^ra_1,\dots,\d^ra_n$, $r>0$. 
An extension $L/K$ of $\d$-fields said to be \emph{$\d$-finitely generated}, if $L$ is the quotient field
of some $\dK$-finitely generated $\dK$-subalgebra of $L$. 

An affine $\d$-group $\G=\op{Spec}_{\d}(H)$ is said to be \emph{$\d$-algebraic} if the $\dK$-Hopf algebra
$H$ is $\d$-finitely generated as a $\dK$-algebra; see Remark \ref{rem:algebraic}. Obviously, ``algebraic" implies ``$\d$-algebraic". 

\begin{lemma}\label{lem:split_by}
Every right (or left) $\dK$-torsor for an affine $\dK$-group $\G$ is split by some $\dK$-field. 
It is split by a $\d$-finitely generated extension $L/K$ of $\d$-fields, 
if $\G$ is $\d$-algebraic.
\end{lemma}
\pf Suppose that $B$ is a right $H$-Galois $\dK$-algebra, as above.
Choose arbitrarily a maximal $\d$-stable ideal $\mathfrak{m}$ of $B$, and construct
$R=B/\mathfrak{m}$, a simple $\dK$-ring. Since $\op{char} K=0$, $R$ is an integral domain by 
\cite[Lemma 1.17]{VS}. The quotient field $L=Q(R)$ of $R$ uniquely turns into a $\dK$-field.
By applying $L \otimes_B$ to $\tilde{\rho}$, it follows that $B$ is split by $L$, proving the first assertion.
If $H$ is $\d$-finitely generated, then $B$ and $R$ are so. It follows that the $L/K$ above is
$\d$-finitely generated, proving the second. 
\epf

\begin{prop}\label{prop:form_of_fd_object}
Suppose that $A$ is $\dK$-object of finite $K$-dimension.
Then $\ms{Aut}_{\d}(A)$ is an affine $\d$-algebraic $\dK$-group, and 
every twisted form of $A$ is split by some $\d$-finitely generated extension $L/K$ of $\d$-fields. 
\end{prop}
\pf 
We have only to prove that $\ms{Aut}_{\d}(A)$ is an affine $\d$-algebraic $\dK$-group, since the rest then follows from
the preceding Lemma.

Choose a $K$-basis $v_1,\dots,v_n$ of $A$. Let 
\[
F=K[x_{ij}, x_{ij}', x_{ij}'', \dots, x_{ij}^{(r)}, \dots ]
\]
denote the free $\dK$-algebra in indeterminates $x_{ij}$, where $1\le i,j\le n$.
Let
\[
G=F_{d}\, (=F[1/d])
\]
denote the localization by the determinant $d=\det X$ of the $n\times n$ matrix $X=\big(x_{ij}\big)_{i,j}$
which has the indeterminates above as entries. 
This $G$ has the $\d$-operator uniquely extending the one $\d x_{ij}^{(r)}=x_{ij}^{(r+1)}$,\ $r\ge 0$, on $F$. 
We have a $G$-linear bijection $\phi: A\ot_K G\to A \ot_K G$ determined by
\[ \phi(v_j \ot 1)= \sum_{i=1}^n v_i\ot x_{ij},\quad 1 \le j \le n. \]
This is alternatively expressed as 
\[
\phi(v_1\ot 1,\dots,v_n\ot 1)=(v_1,\dots,v_n)\ot X 
\]
by matrix presentation;
such presentation will be used in \eqref{eq:matrix_pres1}, \eqref{eq:matrix_pres2} and \eqref{eq:matrix_pres3}, as well. 

Let
\[
H=G/\mathfrak{a},
\]
where $\mathfrak{a}$ is the smallest $\d$-stable ideal of $G$ such that 
the base extension $\phi_H : A\ot_K H \to A \ot_K H$ of $\phi$ along $G \to G/\mathfrak{a}=H$ 
is an endomorphism of the $\d$-$H$-object $A\ot_K H$; obviously, it is necessarily an automorphism.
This $\mathfrak{a}$ is, in fact, given by the relations which ensure that $\phi_H$ commutes with the
structure maps of $A$ (see \cite[Section 7.6]{W}), and with the $\d$-operator. Explicitly, the latter relation for commuting with $\d$-operator is
\begin{equation}\label{eq:matrix_pres1}
XD=DX +\d X,
\end{equation}
where $D \in M_n(K)$ is the matrix determined by
\begin{equation}\label{eq:matrix_pres2}
\d (v_1,\dots,v_n)= (v_1, \dots, v_n )D. 
\end{equation}

We see that $H$ represents the functor $\ms{Aut}_{\d}(A)$ regarded to be set-valued. In fact, 
for every $R \in \dK\text{-}\mathtt{Algebras}$, we have the natural bijection
\[
\op{Spec}_{\d}(H)(R) \to \mr{Aut}_{\d\text{-}R}(A\otimes_K R),\quad f\mapsto~\text{the base extension of}~\phi_H~\text{along}~f. 
\]
By Yoneda's Lemma, $H$ uniquely turns into a $\dK$-Hopf algebra with respect to the familiar Hopf-algebra structure
\begin{equation}\label{eq:matrix_pres3}
\Delta X=X \ot X,\quad \epsilon X=I,\quad \mathcal{S}X=X^{-1},  
\end{equation}
where $\Delta,\ \epsilon$ and $\mathcal{S}$ denote the coproduct, the counit and the antipode, respectively,
and it represents the group-valued functor $\ms{Aut}_{\d}(A)$.
Since $H$ is obviously $\d$-finitely generated, the desired result follows. 
\epf

For $K$ as above, we choose and fix an extension $\mathcal{U}/K$ of $\d$-fields
into which every $\d$-finitely generated extension $L/K$ of $\d$-fields
can be embedded. There exists such an extension; a universal extension \cite[Chapter~III, Section~7]{K1} over $K$
is an example. 

For an affine $\d$-algebraic $\dK$-group $\G$, 
we define $H^1_{\d}(K,\G)$ by
\begin{equation}\label{eq:H1KG}
H^1_{\d}(K,\G):= H^1_{\d}(\mathcal{U}/K, \G). 
\end{equation}

The $\d$-automorphism-group functor $\ms{Aut}_{\d}(\G) : T \mapsto \op{Aut}_{\d\text{-}T}(\G_T)$ 
of the right $\G$-equivariant $\dK$-scheme $\G$ is naturally isomorphic to $\G$ itself;\ $\op{Aut}_{\d\text{-}T}(\G_T)$
consists of the natural automorphisms of the functor $\G_T: \d\text{-}T\text{-}\mathtt{Algebras}\to \mathtt{Groups}$. 
This fact, combined with Lemma \ref{lem:split_by}, shows that 
$H^1_{\d}(K,\G)$ classifies all right $\dK$-torsors for $\G$.

For a $\dK$-object $A$ of finite $K$-dimension, we define
\[ H^1_{\d}(K,\mathsf{Aut}_{\d}(A)):= H^1_{\d}(\mathcal{U}/K, \mathsf{Aut}_{\d}(A)). \]
This classifies all twisted forms of $A$, as is seen from Proposition \ref{prop:form_of_fd_object}.


\subsection{$\dK$-Bi-torsors and bi-Galois $\dK$-algebras}\label{subsec:biGalois}
Let $\G=\op{Spec}_{\d}(H)$ be an affine $\dK$-group. Suppose that $\X=\op{Spec}_{\d}(B)$
is a right $\dK$-torsor for $\G$, or in other words, $B=(B,\rho)$ is a right $H$-Galois $\dK$-algebra. 
Tracing the argument of \cite{Sc} modified into our differential situation, we see that there exists uniquely (up to isomorphism)
a pair $(H', \lambda)$ of a $\dK$-Hopf algebra $H'$ and a left $H'$-comodule $\dK$-algebra structure
$\lambda : B \to H'\ot_K B$ such that (i)~$(B, \lambda)$ is a left $H'$-Galois $\dK$-algebra, and (ii)~$\lambda$ and $\rho$
commute in the sense that 
\begin{equation}\label{eq:lambda_rho}
(\lambda \ot \op{id}_H)\circ \rho = (\op{id}_{H'}\ot \rho)\circ \lambda. 
\end{equation}
We say that $B$ is an \emph{$(H', H)$-bi-Galois $\dK$-algebra}. 
Accordingly, we have uniquely a pair of an affine $\dK$-group $\G'$ and its action on $\X$ from the left, such that
(i)~$\X$ is a left $\dK$-torsor for $\G'$, and (ii)~the actions on $\X$ by $\G'$ and by $\G$ commute with each other. We say that
$\X$ is a \emph{$\dK$-bi-torsor}. 
We write
\begin{equation}\label{eq:HB_GX}
H^B,\quad \G^{\X}
\end{equation}
for $H'$, $\G'$, respectively. If $B$ (or equivalently, $\X$) is trivial, or namely if $B=H$ (or $\X=\G$), then 
$H^B=H$ and $\G^{\X}=\G$. This, applied after base extension to $B$, shows the following; see the
proof of Proposition \ref{prop:interpret} below for detailed argument.

\begin{prop}\label{prop:HB_GX}
$H^B$ and $\G^{\X}$ are $B/K$-forms of $H$ and of $\G$, respectively. 
\end{prop}

With $K$ replaced by a non-zero $\d$-ring $R$, 
the results above remain true if the relevant
$\d$-$R$-Hopf algebra is flat over $R$. We remark that
$\d$-$R$-torsors are then required, in addition to the $\tilde{\rho}$ being isomorphic, to be
faithfully flat over $R$. 


\subsection{Interpretation of $H^1_{\d}(K, \G)\to H^1_{\d}(K, \mathsf{Aut}_{\d}(\g))$}\label{subsec:interpret}
Let $\G=\op{Spec}_{\d}(H)$ be an affine algebraic
$\dK$-group, and set $\g:=\op{Lie}(\G)$. Then $\g$ is a $\dK$-Lie algebra of finite
$K$-dimension, whence the $\d$-automorphism-group functor $\ms{Aut}_{\d}(\g)$ is an affine $\d$-algebraic 
$\dK$-group by Proposition \ref{prop:form_of_fd_object}.
We see that the left adjoint action by $\G$ on $\g$ gives rise to a morphism
of affine $\d$-algebraic $\dK$-groups
\[ \mathsf{Ad} : \G \to \ms{Aut}_{\d}(\g), \]
which induces naturally a map between the cohomology sets
\begin{equation}\label{eq:Ad}
\mathsf{Ad}_* : H^1_{\d}(K, \G) \to H^1_{\d}(K, \ms{Aut}_{\d}(\g)).
\end{equation}

Given a right $\dK$-torsor $\X$ for $\G$, we define
\begin{equation}\label{eq:gX}
\g^{\X}:=\Lie(\G^{\X}). 
\end{equation}
This is a twisted form of $\g=\op{Lie}(\G)$,
since $\G^{\X}$ is a twisted form of $\G$; see Proposition \ref{prop:HB_GX}.

\begin{prop}\label{prop:interpret} 
$\mathsf{Ad}_*$ is interpreted in terms of twisted forms so as
\begin{equation}\label{eq:interpret}
[\text{a right}~\dK\text{-torsor}\ \X \ \text{for}\ \G] \ \mapsto \ [\g^{\X}], 
\end{equation} 
where $[\ \, ]$ indicates isomorphism classes.
\end{prop}
\pf
In this proof we write $\ot$ for $\ot_K$, and use the base-on-left notation
for base extensions; see Remark \ref{rem:right_left}.  

Suppose that $\X =\op{Spec}_{\d}(B)$ is a right $\dK$-torsor for $\G$, or in other words, $B=(B, \rho)$ is
a right $H$-Galois $\dK$-algebra. 

Let $\gamma \in \G(B \ot B)$. This gives the automorphism 
\[
\ell_{\gamma} : (B\ot B)\ot H \overset{\simeq}{\longrightarrow} (B\ot B) \ot H
\] 
of the right $((B\ot B) \ot H)$-Galois $\d\text{-}(B\ot B)$-algebra $(B\ot B)\ot H$ defined by
\[ \ell_{\gamma}((b\ot c) \ot h)=(b\ot c)\gamma(h_{(1)})\ot h_{(2)},\quad b, c \in B,\ h \in H. \]
Here and in what follows, we let
\[ \Delta(h)=h_{(1)}\ot h_{(2)},\quad (\Delta \ot \op{id})\circ \Delta(h)=h_{(1)}\ot h_{(2)}\ot h_{(3)} \]
denote the coproduct on $H$. The \emph{right co-adjoint action} 
\[
\ms{Coad}(\gamma) : (B\ot B)\ot H \overset{\simeq}{\longrightarrow} (B\ot B) \ot H 
\]
by $\gamma$ is defined by
\[ 
\ms{Coad}(\gamma )((b\ot c) \ot h)=(b\ot c)\gamma(h_{(1)})\gamma^{-1}(h_{(3)})\ot h_{(2)}. 
\]
This is an automorphism of the $\d$-$(B\ot B)$-Hopf algebra $(B\ot B)\ot H$. 
Note that $\ell_{\gamma}$ turns into an isomorphism of left $((B\ot B) \ot H)$-Galois $\d\text{-}(B\ot B)$-algebras,
if one twists through $\ms{Coad}(\gamma)$ the obvious co-action by $(B\ot B) \ot H$ on the domain.
Explicitly, this means that
\begin{equation}\label{eqref:ellgamma}
(\ms{Coad}(\gamma )\, \ot_{B\ot B} \, \ell_{\gamma})\circ \Delta_{(B\ot B)\ot H}
= \Delta_{(B\ot B)\ot H} \circ \ell_{\gamma} \ \, \text{on} \ \, (B\ot B)\ot H,
\end{equation}
where $\Delta_{(B\ot B)\ot H}$ denotes the coproduct on $(B\ot B)\ot H$. 

Suppose that the $\gamma$ above is a cocycle for computing $H^1_{\d}(B/K,\G)$ which gives the $B/K$-form
$B$ through $\tilde{\rho}$. 
This means that the commutative diagram 
\[
\begin{xy}
(0,0)   *++{(B\ot B) \ot B}  ="1",
(-23,-16)  *++{(B\ot B) \ot H}    ="2",
(23,-16)   *++{(B\ot B) \ot H}  ="3",
{"1" \SelectTips{cm}{} \ar @{->}_{d_1\tilde{\rho}} "2"},
{"1" \SelectTips{cm}{} \ar @{->}^{d_2\tilde{\rho}}"3"},
{"2" \SelectTips{cm}{} \ar @{->}_{\ell_{\gamma}} "3"}
\end{xy}
\]
of right $((B\ot B) \ot H)$-Galois $\d\text{-}(B\ot B)$-algebras, where $d_i$, $i=1,2$, denote
the base extensions along
\[ B\to B\ot B,\quad b\mapsto 1\ot b,\ b\ot 1. \]

Recall that $B$ is an $(H^B,H)$-bi-Galois $\d$-$K$-algebra. 
By \cite[Theorem 3.5]{Sc}, 
the Hopf algebra $H^B$ consists of the elements $\sum_i b_i \ot c_i$ in $B\ot B$ such that
\begin{equation}\label{eq:coinv}
\sum_i (b_i)_{(0)}\ot (c_i)_{(0)}\ot (b_i)_{(1)}(c_i)_{(1)}=\sum_i b_i\ot c_i \ot 1\ \, \text{in} \ \, (B \ot B) \ot H, 
\end{equation}
where $\rho(b)=b_{(0)}\ot b_{(1)}$. Moreover,
\begin{equation}\label{eq:mu}
\mu : B \ot H^B \to B\ot B,\quad \mu(b \ot z)= bz 
\end{equation}
is an isomorphism of left $(B \ot H^B)$-Galois $\d$-$B$-algebras. 
Define 
\[ 
\nu :=\tilde{\rho}\circ \mu : B\ot H^B \to B\ot H. 
\]

Recall from Section \ref{subsec:biGalois} uniqueness of the pair $(H', \lambda)$, and apply it first over $B$, and next over $B \ot B$. 
Then one sees the following. First, there uniquely exists an isomorphism $\theta: B\ot H^B \to B\ot H$ of $\d$-$B$-Hopf algebras
such that
\[ (\theta \ot \nu)\circ \Delta_{B\ot H^B}= \Delta_{B\ot H}\circ \nu, \]
where $\Delta_{B\ot H^B}$ and $\Delta_{B\ot H}$ denote the coproducts on the $\d$-$B$-Hopf algebras. 
In fact, this $\theta$ is the unique isomorphism between the two $\d$-$B$-Hopf-algebras that is compatible with their
co-actions on $B \ot B$. (Notice that this $\theta$ ensures Proposition \ref{prop:HB_GX}.)
Next, the last commutative diagram, with $((B\ot B) \ot H)^{(-)}$ applied (see \eqref{eq:HB_GX}), induces the commutative diagram
\[
\begin{xy}
(0,0)   *++{(B\ot B) \ot H^B}  ="1",
(-23,-16)  *++{(B\ot B) \ot H}    ="2",
(23,-16)   *++{(B\ot B) \ot H}  ="3",
{"1" \SelectTips{cm}{} \ar @{->}_{d_1\theta} "2"},
{"1" \SelectTips{cm}{} \ar @{->}^{d_2\theta}"3"},
{"2" \SelectTips{cm}{} \ar @{->}_{\ms{Coad}(\gamma )} "3"}
\end{xy}
\]
of $\d\text{-}(B\ot B)$-Hopf algebras; notice from \eqref{eqref:ellgamma} that 
$\ell_{\gamma}$
induces $\ms{Coad}(\gamma)$. 

Notice from \eqref{eq:LieG} that $\g^{\X}=\op{Der}_{\epsilon}(H^B,K)$. 
Then one sees that $\theta$ induces an isomorphism 
\[ \theta^* : B \ot \g \overset{\simeq}{\longrightarrow} B \ot \g^{\X} \]
of $\d$-$B$-Lie algebras. Moreover, the last commutative diagram induces by duality the commutative diagram
\[
\begin{xy}
(0,0)   *++{(B\ot B) \ot \g^{\X}}  ="1",
(-23,-16)  *++{(B\ot B) \ot \g}    ="2",
(23,-16)   *++{(B\ot B) \ot \g}  ="3",
{"1" \SelectTips{cm}{} \ar @{<-}_{d_1(\theta^*)} "2"},
{"1" \SelectTips{cm}{} \ar @{<-}^{d_2(\theta^*)}"3"},
{"2" \SelectTips{cm}{} \ar @{<-}_{\ms{Ad}(\gamma^{-1})} "3"}
\end{xy}
\]
of $\d\text{-}(B\ot B)$-Lie algebras, where the horizontal arrow indicates the left adjoint action by 
$\gamma^{-1}$. 
We may reverse the direction of the arrow, changing the label into the left adjoint action $\ms{Ad}(\gamma)$ by $\gamma$. 
The result shows that $\ms{Ad}(\gamma)$, regarded as
a cocycle for computing $H^1_{\d}(B/K, \ms{Aut}_{\d}(\g))$, 
gives the twisted form $\g^{\X}$ of $\g$ which is split by $B$, indeed.

Recall from \eqref{eq:H1KG} the definition $H^1(K,\G):=H^1(\mathcal{U}/K,\G)$.
Let $\psi$ is an element of $H^1(K,\G)$. Then this arises from a cocycle
$\gamma$ such as above, 
which gives a $\d$-$K$-torsor $\X=\op{Spec}_{\d}(B)$ for $\G$, through 
a $\d$-$K$-algebra map $j: K \to \mathcal{U}$. Thus, $\psi$ is represented
by the cocycle given as the composite
\begin{equation*}\
H \overset{\gamma}{\longrightarrow} B\otimes B\overset{j\otimes j}{\longrightarrow} 
\mathcal{U}\otimes \mathcal{U}.
\end{equation*}
This cocycle is seen to give
the $\mathcal{U}/K$-form $B$ of the right $H$-comodule $\d$-$K$-algebra $H$. The argument in the preceding paragraphs shows
that $\ms{Ad}_*(\psi)$ is represented by the base extension of the automorphism 
$\ms{Ad}(\gamma)$ along the $\d$-$K$-algebra map $j\otimes j$.
This base extension is seen to be a cocycle which gives the $\mathcal{U}/K$-form 
$\g^{\X}$ of $\g$. This completes the proof.
\epf

Let $\G=\op{Spec}_{\d}(H)$ be an affine algebraic $\dK$-group with $\g=\Lie(\G)$, as above.
Recall from \eqref{eq:LieG} that $\g=\op{Der}_{\epsilon}(H, K)$. 

\begin{prop}\label{prop:dadD}
Regard $H$ merely as the trivial right $H$-Galois $K$-algebra, forgetting $\d$ on it. 
\begin{itemize}
\item[\textup{(1)}]
Given an element $D \in \g$, define 
\begin{equation}\label{eq:dD}
\d_D : H \to H,\quad \d_D(h)=\d h+ D(h_{(1)})h_{(2)},  
\end{equation}
where $\d$ denotes the original operator on $H$. Then this is a $\d$-operator 
with which $H$ is made into a right $H$-Galois $\dK$-algebra. Conversely, such a $\d$-operator uniquely arises in this way.
\item[\textup{(2)}]
Given  an element $D\in \g$, let $\X_D$ denote the right $\dK$-torsor for $\G$ which is represented by the right $H$-Galois $\dK$-algebra
$(H, \d_D)$ obtained above. Then the twisted form $\g^{\X_D}$ of $\g$ is the $K$-Lie algebra $\g$ equipped with the new $\d$-operator
\[ \d + \ms{ad}(D) : \g\to \g,\quad  z \mapsto \d z +[D, z], \]
where $\d$ denotes the original operator on $\g$. Thus $\g^{\X_D}$ is quasi-isomorphic to the original $\g$;
see Definition \ref{def:quasi-isom}. 
\end{itemize}
\end{prop}
\pf
(1)\
Suppose that $\d_1$ is a desired operator, or namely, $(H,\d_1)$ is a right $H$-Galois $\dK$-algebra.
Then one sees that
$\d_1-\d : H \to H$ is a $K$-linear derivation and is at the same time a right $H$-comodule map. 
It follows that $\d_1$ is necessarily of
the form $\d_D$ with $D\in \g$ uniquely determined. 
Such $\d_D$ is seen to be a desired operator for any $D$, indeed. 

(2)\ 
Let $H'=H^{(H, \d_D)}$. Then $\g^{\X_D}=\op{Der}_{\epsilon}(H', K)$. 
Using the uniqueness of the $\dK$-Hopf algebra $H'$ in general, which was discussed in Section \ref{subsec:biGalois}, 
we see that the present $H'$ is the $K$-Hopf algebra $H$ equipped with the $\d$-operator
\[ H \to H,\quad h\mapsto \d h+ D(h_{(1)})h_{(2)}- h_{(1)}D(h_{(2)}). \]
This implies the desired result. 
\epf

A simple consequence of the proposition above is the following.

\begin{corollary}\label{cor:trivial_form}
Let $\g$ be a $\dK$-Lie algebra of finite $K$-dimension. 
Once the Lie algebra $\Lie(\G)$ of some 
affine algebraic $\dK$-group $\G$ is shown to be a twisted form of $\g$, then every $\dK$-Lie algebra quasi-isomorphic to 
$\Lie(\G)$ is a twisted form of $\g$, as well. 
\end{corollary}


\subsection{Differential $\dK$-objects arising from $C$-linear objects}
Let $K$ be a $\d$-field of characteristic zero. Let
\[ C=C_K\, (=\{ \, x\in K \mid \d x=0\, \}) \]
denote the field of constants in $K$, which is necessarily of characteristic zero.  
In this subsection we let $\ot$ denote the tensor product $\ot_C$ over $C$.

Let $A_0$ be a $C$-linear object. We can and do regard the base extension $A_0\ot K$ as
a $\dK$-object with respect to the operator $\d_0$ defined by
\begin{equation*}\label{eq:d0}
\d_0 : A_0\ot K \to A_0\ot K,\quad a \ot x \mapsto a \ot \d x.
\end{equation*}
For every $\dK$-algebra $R$, $A_0 \ot R$ is similarly a $\d$-$R$-object, and is a base extension of 
the $\dK$-object $A_0\ot K$ above. 

\begin{prop}\label{prop:Aut}
If the automorphism-group functor $\ms{Aut}(A_0)$ of $A_0$ happens to be an affine $C$-group,
represented by a $C$-Hopf algebra $H_0$, then  
the $\d$-automorphism-group functor $\ms{Aut}_{\d}(A_0\ot K)$ of the $\dK$-object
$(A_0\ot K, \d_0)$ is an affine $\dK$-group, represented by the $\dK$-Hopf algebra $H_0\ot K$. 
\end{prop}
\pf
Let $R\in \dK\text{-}\mathtt{Algebras}$. 
One sees that every automorphism of the $\d$-$R$-object $A_0 \ot R$ restricts to an automorphism of $A_0 \ot C_R$ over 
the $C$-algebra $C_R$ of constants in $R$, and so it is uniquely 
presented as the base extension of the restriction.
This shows $\mathrm{Aut}_{\d\text{-}R}(A_0\ot R)= \mathrm{Aut}_{C_R}(A_0\ot C_R)$; this last is naturally isomorphic to 
$\op{Spec}_C(H_0)(C_R)=\op{Spec}_{\dK}(H_0\ot K)(R)$. This proves the proposition. 
\epf

We remark that the proposition follows from the proof of Proposition \ref{prop:form_of_fd_object} if 
$(\dim_K(A_0\ot K)=)\dim_CA_0<\infty$. 
For the relation \eqref{eq:matrix_pres1} turns into $\d X=O$ since $D=O$. 

\begin{corollary}
If $\dim_C A_0<\infty$, then $\ms{Aut}_{\d}(A_0\ot K)$ is an affine algebraic $\dK$-group.
\end{corollary}
\pf 
This follows from the proposition above, since $\ms{Aut}(A_0)$ is an affine algebraic $C$-group under the
assumption; see \cite[Section 7.6]{W}. 
\epf

The following result would be worth presenting, though it is not essentially used in this paper. 

\begin{prop}
Assume $\dim_C A_0<\infty$, and that the field $C$ is algebraically closed. 
Then every twisted form of the $\dK$-object $(A_0 \ot K, \d_0)$ is split
by some (finitely generated) Picard-Vessiot extension $L$ over $K$. 
\end{prop}
\pf
By the preceding results the first assumption implies
$\ms{Aut}(A_0)=\op{Spec}_C(H_0)$ and $\ms{Aut}_{\d}(A_0\ot K)=\op{Spec}_{\dK}(H_0\ot K)$,
where $H_0$ is a finitely generated $C$-Hopf algebra. 

Let $B$ be a right $(H_0\ot K)$-Galois $\dK$-algebra, and regard it 
as a twisted form of the right $(H_0\ot K)$-comodule $\dK$-algebra $H_0\ot K$. 
We should prove that this twisted form $B$ is split by some 
$L/K$ as above.
It suffices to prove that there exists a $\dK$-algebra map from $B$ to such an $L$, since
$B$ is split by $B$, itself. 

We have the $\d$-$B$-algebra isomorphism 
\[
\tilde{\rho} : B\ot_K B \overset{\simeq}{\longrightarrow} B\ot_K(H_0\ot K)=B\ot H_0 
\]
as in \eqref{eq:rho}. 
Choose a simple quotient $\dK$-algebra $R$ of $B$, as in the proof of Lemma \ref{lem:split_by}. 
Then $R$ is an integral domain by \cite[Lemma 1.17]{VS}, as before. 
This is finitely generated 
as a $K$-algebra since $B$ is such. 
The quotient field $L=Q(R)$ of $R$ uniquely turns into a $\d$-field, which is necessarily a finitely generated extension over $K$. 
The second assumption above,  combined with \cite[Lemma 4.2]{AMT}, implies that the field $C_L$ of constants in $L$ equals $C$. 
This $L/K$ will be proved to be a desired Picard-Vessiot extension by \cite[Definition 1.8 and Theorem 3.11]{AMT}, 
if one sees that the canonical $\d$-$R$-algebra map $R \ot C_{R\ot_K R} \to R\ot_K R$ which arises from 
the embedding $C_{R\ot_K R} \hookrightarrow R\ot_K R$ of the constants into $R\ot_K R$ is surjective. 
(By \cite[Proposition 6.7]{AMT} this canonical map is injective, though this fact is not needed here.)
Indeed, the desired surjectivity is seen from the commutative diagram
\[
\begin{xy}
(0,16)   *++{B\ot H_0}  ="1",
(30,16)  *++{B\ot_K B}    ="2",
(0,0)   *++{R\ot C_{R\ot_KR}}  ="3",
(30,0)   *++{R\ot_KR.}  ="4",
{"1" \SelectTips{cm}{} \ar @{->}_{\simeq}^{\tilde{\rho}^{-1}} "2"},
{"1" \SelectTips{cm}{} \ar @{->} "3"},
{"2" \SelectTips{cm}{} \ar @{->>}^{} "4"},
{"3" \SelectTips{cm}{} \ar @{->}^{} "4"}
\end{xy}
\]
Here the vertical arrow on the left-hand side naturally arises from the composite of $\tilde{\rho}^{-1}|_{H_0}:H_0\to B \ot_K B$
with the natural surjection $B\ot_KB\to R\ot_KR$, which clearly takes values in $C_{R\ot_K R}$. 
\epf


\section{Proof of the theorem and computations}\label{sec4}
Throughout in this section we let $K:=\C$, and write $\ot$ for $\ot_K$. 

Suppose that we are in the situation of Section \ref{sec:intro}.
Let $\g_0$ be a complex simple Lie algebra, and let $\g=\g_0(K)$ denote the $\dK$-Lie algebra
as in \eqref{eq:Lie_functor}. 

\subsection{Two key facts}
One key fact for us is the following description of the automorphism-group scheme $\ms{Aut}(\g_0)$ of $\g_0$. 
Recall that the finite group $\Gamma=\ms{Out}(\g_0)$ of outer-automorphisms of $\g_0$
is explicitly given by \eqref{eq:Gamma}; this $\Gamma$ will be identified with the associated, finite constant group scheme. 
Let $\G_0^{\circ}$ be the adjoint simple $\c$-group 
associated with $\g_0$. A natural action by $\Gamma$ on $\G_0^{\circ}$ constitutes 
an affine algebraic $\c$-group
\[ \G_0 = \G_0^{\circ}\rtimes \Gamma \]
of semi-direct product, so that
\[ \Lie(\G_0)=\Lie(\G_0^{\circ})=\g_0, \]
and the adjoint action by $\G_0$ on $\g_0$ gives an isomorphism
\begin{equation}\label{eq:ad_iso}
\ms{Ad} : \G_0 \overset{\simeq}{\longrightarrow} \ms{Aut}(\g_0)
\end{equation}
of affine algebraic $\c$-groups.  
By restriction this $\ms{Ad}$ induces the identity $\Gamma = \ms{Out}(\g_0)$. 
Note that $\G_0^{\circ}$ is the connected component of $\G_0$ containing the identity element. 
See \cite[Chapter 4, Section 4,\ $1^{\circ}$]{OV}.

Suppose $\G_0=\op{Spec}_{\c}(H_0)$, and define 
\[ \G=\op{Spec}_{\d\text{-}K}(H_0\otimes_{\c} K). \]
Then one sees $\g=\Lie(\G)$. Moreover, it follows from \eqref{eq:ad_iso} and Proposition \ref{prop:Aut} that
the adjoint action by $\G$ on $\g$ gives an isomorphism
\[
\ms{Ad} : \G \overset{\simeq}{\longrightarrow} \ms{Aut}_{\d}(\g)
\]
of affine algebraic $\d$-$K$-groups.  
This together with Proposition \ref{prop:interpret} prove the following. 

\begin{prop}\label{prop:describe}
Every form of the $\dK$-Lie algebra $\g$ uniquely arises, as described by \eqref{eq:interpret},
from a right $\dK$-torsor for $\G$.
\end{prop}

Another key fact is the cohomology vanishing of the (non-differential) 
Amitsur 1st cohomology due to Steinberg (see Serre \cite[III, 2.3, Theorem 1$'$]{Se}), 
\begin{equation}\label{eq:vanishing}
H^1(K, \ms{F})=0,
\end{equation}
where $\ms{F}$ is a connected affine algebraic $K$-group. 
This is proved more generally when $K$ is replaced by a perfect field, say $K'$, of dimension $\le 1$ 
\cite[Definitoion on Page 78]{Se}, and in addition, $\ms{F}$ is assumed to be smooth if $\op{char} K'>0$;
note that every affine algebraic $K$-group is necessarily smooth since $\op{char} K=0$.  
One sees that $K\, (=\C)$ is a $\mathrm{(C_1)}$-field by Tsen's Theorem, whence $K$ is of dimension $\le 1$ by \cite[Corollary on Page 80]{Se}. 

\subsection{Proof of Theorem \ref{thm1}, Part 1: Case $\Gamma=\{1\}$} 
In this case, $\G$, regarded as an affine $K$-group, is connected. 
By \eqref{eq:vanishing} applied to this $\G$, it follows that
every right $\dK$-torsor for $\G$, regarded as a right $K$-torsor, is trivial. 
Propositions \ref{prop:dadD} and  \ref{prop:describe} conclude the proof.

\subsection{Galois descent}\label{subsec:Galois_descent}
To proceed to Parts 2 and 3, suppose that we are now in Case $\Gamma \ne \{1\}$. 

Note that $\Gamma$, regarded as a finite constant $\c$-group scheme,
is represented by the dual $(\c \Gamma)^*$
of the group algebra $\c \Gamma$; this $(\c \Gamma)^*$ is the separable part $\pi_0(H_0)$ of $H_0$ \cite[Page 49]{W}, that is, the largest separable
subalgebra (in fact, Hopf subalgebra) of the $\c$-Hopf algebra $H_0$. 
Suppose $\G_0^{\circ}=\op{Spec}_{\c}(J_0)$, and define
\[ H:=H_0\ot_{\c} K,\quad J:=J_0\ot_{\c}K,\quad Z:=(\c \Gamma)^* \ot_{\c} K\, (=(K\Gamma)^*), \]
which are naturally $\dK$-Hopf algebras, such that $\G=\op{Spec}_{\d}(H)$, in particular.  
One sees that $Z\subset H$ is a $\dK$-Hopf subalgebra, and
\begin{equation}\label{eq:J}
J=H/(Z^+),
\end{equation}
where $(Z^+)$ is the ideal (in fact, $\d$-stable Hopf ideal) generated by the augmentation ideal $Z^+=\op{Ker}(\epsilon:Z\to K)$
of $Z$, that is, the kernel of the counit. 
Since $\Gamma$ acts innerly on $\G_0^{\circ}\, (\subset \G_0)$ from the right, it acts from the left
on $J$ as $\dK$-Hopf-algebra automorphisms. The action gives rise by adjoint to the co-action
$J \to J \ot Z$ by $Z=(K\Gamma)^*$, so that the associated smash coproduct $Z\rcosmash J$ (see \cite[Definition 10.6.1]{Mo})
coincides with $H$.
Here one should recall $\G_0=\Gamma \ltimes \G_0^{\circ}\, (=\G_0^{\circ}\rtimes \Gamma)$. 

Choose arbitrarily a right $\dK$-torsor $\X=\op{Spec}_{\d}(B)$ for $\G=\op{Spec}_{\d}(H)$. In view of
Proposition \ref{prop:describe} we wish to describe the $\dK$-Lie algebra $\g^{\X}\, (=\Lie(\G^{\X}))$. 
Let $H':=H^B$, or in other words, suppose $\G^{\X}=\op{Spec}_{\d}(H')$, so that
$B$ is an $(H',H)$-bi-Galois $\dK$-algebra. We are going to prove the following.

\begin{prop}\label{prop:key}
There is a finite-dimensional $\dK$-Hopf subalgebra $Z'$ of $H'$ such that 
\begin{itemize}
\item[\textup{(i)}] $Z'$ is separable as a $K$-algebra; 
\item[\textup{(ii)}] the associated quotient $\dK$-Hopf algebra 
\begin{equation}\label{eq:J'}
J':=H'/(Z'^+)\quad (\text{cf}.~\eqref{eq:J}) 
\end{equation}
has the trivial separable part, $\pi_0(J')=K$, or in other words, it includes no non-trivial
separable $K$-subalgebra. 
\end{itemize}
\end{prop}

This implies that the affine $K$-group $\op{Spec}(H')$ includes $\op{Spec}(J')$ as the connected component
containing the identity element, and thereby concludes
\begin{equation}\label{eq:conclude}
\g^{\X}=~\text{the Lie algebra of the affine}~\dK\text{-group}~\op{Spec}_{\d}(J')
\end{equation}
as $\dK$-Lie algebras. Therefore, we aim first to prove the proposition above, and then to describe the $\g^{\X}$ above. 

Let $\rho : B \to B \ot H$ denote the structure map on $B$, and define
\[ R := \rho^{-1}(B \ot Z). \]
Then this $R$ is a right $\dK$-Galois algebra for $Z$, or in other words, $\op{Spec}_{\d}(R)$ is a right $\dK$-torsor
for the finite constant $\dK$-group scheme $\Gamma_K$ given by $\Gamma$; it arises from the
the right $\dK$-torsor $\X=\op{Spec}_{\d}(B)$ for $\G$ through the restriction map 
$H^1_{\d}(K, \G)\to H^1_{\d}(K,\Gamma_K)$ which is defined in the differential situation, as well, just as in the ordinary situation. 
Note that $R$ is naturally a $\dK$-algebra of finite $K$-dimension, and is 
a \emph{Galois $K$-algebra} with Galois group $\Gamma$ in the classical sense that
the $K$-algebra map $R \rtimes \Gamma \to \op{End}_{K}(R)$ which arises from the natural module-action on $R$
by the semi-direct product $R\rtimes \Gamma$ is an isomorphism. 
Note that $R\rtimes \Gamma$ is naturally a non-commutative $\d$-$K$-algebra with $\Gamma\, (=\{1\}\times \Gamma)$ 
included in constants.  
A $\d$-$(R\rtimes \Gamma)$-\emph{module} is thus an $R$-module $M$ equipped with 
an additive operator $\d$ and a $\Gamma$-action of $K$-linear automorphisms, such that
\[
\d(\gamma m)=\gamma(\d m),\quad \d(am)=(\d a)m+a(\d m),\quad \gamma(am)=(\gamma a)(\gamma m), 
\]
where $\gamma \in \Gamma$, $a \in R$ and $m \in M$. We call this a \emph{$(\d,\Gamma)$-$R$-module},
to treat $\d$ and $\Gamma$ on an equality, 
and let $(\d,\Gamma)\text{-}R\text{-}\mathsf{Modules}$ denote the category of those modules.
The classical Galois Descent Theorem (see \cite[Section 18]{KMRT}) tells us that the
functor $M \mapsto M^{\Gamma}$,\ $\Gamma$-invariants in $M$, 
gives the category equivalence
\[
(\d,\Gamma)\text{-}R\text{-}\mathsf{Modules}\overset{\approx}{\longrightarrow}
\d\text{-}K\text{-}\mathsf{Modules},
\]
whose quasi-inverse is given by the base-extension functor $\ot_KR$. 
In fact, this is a symmetric tensor equivalence, so that there is induced the category equivalence
between their (commutative-)algebra objects, or between any
other kind of linear objects. The category on the left-hand side has the tensor product $\otimes_R$, the unit object $R$ and the
obvious symmetry, while the category on the right-hand side has the tensor product $\ot_K$, the unit object $K$ and
the obvious symmetry. A commutative algebra in $(\d,\Gamma)\text{-}R\text{-}\mathsf{Modules}$ will be called
a $(\d,\Gamma)\text{-}R$-\emph{algebra}; it descends to a $\dK$-algebra by the category equivalence above. 
Similarly, a $(\d,\Gamma)\text{-}R$-\emph{Hopf} or \emph{Lie algebra} is defined, and it descends to a $\dK$-object. 


We have the commutative diagram 
\[
\begin{xy}
(0,16)   *++{B\ot B}  ="1",
(30,16)  *++{B\ot H}    ="2",
(0,0)   *++{R\ot R}  ="3",
(30,0)   *++{R\ot Z}  ="4",
(15,-16) *++{R}="5",
{"1" \SelectTips{cm}{} \ar @{->}_{\simeq}^{\tilde{\rho}} "2"},
{"3" \SelectTips{cm}{} \ar @{^(->} "1"},
{"4" \SelectTips{cm}{} \ar @{^(->}^{} "2"},
{"3" \SelectTips{cm}{} \ar @{->}_{\simeq} "4"},
{"3" \SelectTips{cm}{} \ar @{->}_{\text{mult}} "5"},
{"4" \SelectTips{cm}{} \ar @{->}^{\text{id}_R \ot \, \epsilon} "5"}
\end{xy}
\]
of $\dK$-algebras, where
the upper horizontal arrow indicates the isomorphism $\tilde{\rho}$ (see \eqref{eq:rho}) 
associated with the structure 
map $\rho : B \to B \ot H$ on $B$, and 
the lower one is the analogous isomorphism for the right $Z$-Galois $\dK$-algebra $R$. In addition, $\text{mult} : R\ot R \to R$
indicates the multiplication $x \ot y\mapsto xy$. 
By the base extensions along the two diagonal arrows $\text{mult}: R \ot R \to R$ and $\op{id}_R \ot \, \epsilon: R \ot H \to R$, 
the $\tilde{\rho}$ induces the isomorphism
\begin{equation}\label{eq:BBBJ}
B\ot_R B \overset{\simeq}{\longrightarrow} B\otimes J =B\ot_R(J\ot R). 
\end{equation}

Recall that $\Gamma$ acts on $J$ as $\dK$-Hopf algebra automorphisms. Then one sees that $J\ot R$ is a
$(\d,\Gamma)\text{-}R$-Hopf algebra, 
and hence descends to a $\dK$-Hopf algebra
\[ \mathcal{J}:=(J\ot R)^{\Gamma}. \]
The composite 
\[
B \to B\ot H\to B\ot J=B\ot_R(J\ot R) 
\]
of the structure map on $B$ with
the natural surjection onto $B \ot_R(J\ot R)$
is a $(\d,\Gamma)$-$R$-algebra map, and hence descends to a $\dK$-algebra
map $B^{\Gamma}\to B^{\Gamma}\ot \mathcal{J}$, which we call $\varrho$.

\begin{lemma}\label{lem:BGamma}
$B^{\Gamma}$ is a right $\mathcal{J}$-Galois $\dK$-algebra by the $\varrho$ above. 
\end{lemma}
\pf
One sees that $\varrho$ satisfies the co-associativity and the counit property since 
the last composite does. One sees that \eqref{eq:BBBJ} is an isomorphism of $(\d,\Gamma)$-$R$-algebras, and descends to 
$\tilde{\varrho}: B^{\Gamma}\ot B^{\Gamma}\to B^{\Gamma}\ot \mathcal{J}$, which is, therefore, an isomorphism. 
\epf

Recall $H'=H^B$. Define
$Z':=Z^R$,
so that $R$ is a $(Z',Z)$-bi-Galois $\dK$-algebra. 

\begin{lemma}\label{lem:Z'}
$Z'$ is a finite-dimensional $\dK$-Hopf subalgebra of $H'$ which has the property \textup{(i)}
of Proposition \ref{prop:key}, that is, $Z'$ is separable as a $K$-algebra.  
\end{lemma}
\pf
By \eqref{eq:coinv} we have $Z'\subset H'$. This inclusion is compatible with the Hopf-algebra structure maps,
as is seen from the construction of $H^B$ given in \cite[Theorem 3.5]{Sc}. To verify this here only for the coproduct, 
recall from \eqref{eq:mu} that $H'\subset B\ot B$ gives rise to a left $B$-linear isomorphism $B \ot H'=B\ot B$.
Therefore, we have 
\[
H'\ot H'\subset B\ot H'\ot H'=B\ot B\ot H'=B\ot B\ot B.
\]
The construction cited above tells us that
the coproduct on $H'$ is the restriction of 
\[
B\ot B\to B\ot B\ot B,\quad b\ot c \mapsto b\ot 1\ot c.
\]
This, combined with
the analogous restriction 
of $R\ot R\to R\ot R\ot R$ to the coproduct $Z'\to Z\ot Z'$, shows the desired compatibility, as is verified by 
a commutative diagram in cube. 

The $K$-algebras $Z$, $R$ and $Z'$ turn to be mutually isomorphic
after base extension to some algebraically closed field.
It follows that $Z'$ is finite-dimensional separable, since $Z$ is. 
\epf

Define $J':=H'/(Z'^+)$, as in \eqref{eq:J'}. 
The proof of Proposition \ref{prop:key} completes by proving the next lemma. The following proposition
describes the $\dK$-Lie algebra $\g^{\X}$; see \eqref{eq:gX}.

\begin{lemma}\label{lem:BbiGalois}
$B^{\Gamma}$ is a $(J',\mathcal{J})$-bi-Galois $\dK$-algebra, and $J'$ has the property \textup{(ii)} 
of Proposition \ref{prop:key}, that is, $\pi_0(J')=K$.  
\end{lemma}
\pf
The same argument as proving Lemma \ref{lem:BGamma} shows that $B^{\Gamma}$ is a left $J'$-Galois $\dK$-algebra.
Here one should notice that $\Gamma$ acts (or $Z$ co-acts) trivially on $H'$, and hence on $J'$. 
Indeed, $B^{\Gamma}$ is bi-Galois, since the structure maps 
\[
H'\ot B \leftarrow B \to B\ot H
\]
on $B$ commute with each other
(see \eqref{eq:lambda_rho}), and hence those on $B^{\Gamma}$ do.  

Note that $\pi_0(J)\, (=\pi_0(J_0)\ot_{\c} K)$ equals $K$. 
This is equivalent to saying that the $K$-algebra $J$
contains no non-trivial idempotent even after base extension to some (or any) algebraically closed field. 
It follows that $\mathcal{J}$ and $J'$ have the same property, since $J$ and $\mathcal{J}$, as well as $\mathcal{J}$ and $J'$,
are mutually isomorphic after base extension such as above. 
\epf

Since $\g_0(R)=\g_0\ot_{\c} R$, on which $\Gamma$ acts diagonally, is a $(\d,\Gamma)\text{-}R$-Lie algebra, it descends to
$\g_0(R)^{\Gamma}$, a $\dK$-Lie algebra. Our aim of this subsection is achieved by the following. 

\begin{prop}\label{prop:gX}
The $\dK$-Lie algebra $\g^{\X}$ is quasi-isomorphic to $\g_0(R)^{\Gamma}$.   
\end{prop}
\pf
Recall \eqref{eq:conclude} and the result of Proposition \ref{prop:key} that $B^{\Gamma}$ is a $(J',\mathcal{J})$-bi-Galois $\dK$-algebra. 
By Steinberg's Cohomology-Vanishing \eqref{eq:vanishing}
applied to the connected affine $K$-group $\operatorname{Spec}(\mathcal{J})$, we see that the right $\mathcal{J}$-Galois $K$-algebra 
$B^{\Gamma}$ is isomorphic to $\mathcal{J}$. This together with Proposition \ref{prop:dadD} prove the desired result. 
\epf

We add an important consequence.

\begin{corollary}\label{cor:form}
The twisted forms of $\g_0(K)$ are precisely the $\dK$-Lie algebras quasi-isomorphic to 
$\g_0(R)^{\Gamma}$, where $R$ ranges over all right $(K\Gamma)^*$-Galois $\dK$-algebras. 
\end{corollary}
\pf 
By Propositions \ref{prop:describe} and \ref{prop:gX}, every twisted form is quasi-isomorphic to some $\g_0(R)^{\Gamma}$.
Conversely, any $\g_0(R)^{\Gamma}$ is clearly a twisted form, whence any one that is quasi-isomorphic 
to $\g_0(R)^{\Gamma}$ is, as well, by \eqref{eq:conclude} and Corollary \ref{cor:trivial_form}. 
\epf

Before proceeding we make the following remark: given an integer $n \ge 2$, let $\Lambda_n=K^{\times}/(K^{\times})^n$
denote the quotient group of the multiplicative group $K^{\times}$ by the subgroup of all $n$-th powers. This is an infinite group,
as is easily seen.
Removing the identity element, let $\Lambda_n^{+}=\Lambda_n\setminus \{ 1 \}$. 
Needed here is the set only in $n=2, 3$. The set $\Lambda_2^{+}$ parametrizes the quadratic field extensions over $K$, 
while the set $\Lambda_3^{+}$ modulo the equivalence relation $x\sim x^{\pm 1}$ 
parametrizes
the cubic Galois field extensions over $K$. Therefore, these two classes of field extensions both consist of infinitely many ones.

\subsection{Proof of Theorem \ref{thm1}, Part 2: Case $\Gamma=\mathbb{Z}_2$}
In this case, the right $(K\Gamma)^*$-Galois $\dK$-algebras $R$ are precisely 
\begin{itemize}
\item[(i)] the trivial
one $(K\Gamma)^*$ (equipped with the obvious $\d$-operator), and 
\item[(ii)] the quadratic field extensions over $K$
(equipped with the $\d$-operator uniquely extending the one on $K$).
\end{itemize} 
By Corollary \ref{cor:form}
it remains to show that 
for $R=(K\Gamma)^*$ in (i), we have $\g_0(R)^{\Gamma}\simeq \g_0(K)$. 
Let $\op{Map}(\Gamma, \g_0(K))$ denote the $\Gamma$-set of all maps $\Gamma \to \g_0(K)$, 
equipped with the action 
\[ \gamma f : \gamma' \mapsto f(\gamma'\gamma), \]
where $\gamma, \gamma' \in \Gamma$ and $f \in \op{Map}(\Gamma, \g_0(K))$. 
Regard this naturally as the direct product of $\#\Gamma$-copies of the $\dK$-Lie algebra $\g_0(K)$. 
Then we see that associating to $x \ot a \in \g_0 \ot_{\c} (K\Gamma)^*$, the map $\gamma \mapsto \gamma x \ot a(\gamma)$
gives a $\Gamma$-equivariant isomorphism 
\begin{equation*}
\g_0(R) \overset{\simeq}{\longrightarrow} \op{Map}(\Gamma, \g_0(K)) 
\end{equation*} 
of $\dK$-Lie algebras, whose restriction to the $\Gamma$-invariants is the desired $\g_0(R)^{\Gamma}\simeq \g_0(K)$. 
This completes the proof. 

\subsection{Proof of Theorem \ref{thm1}, Part 3: Case $\Gamma=\mathfrak{S}_3$}
In this case, let $R$ be a right $(K\Gamma)^*$-Galois $\dK$-algebra. 
In view of Corollary \ref{cor:form} we wish to show that $\g_0(R)^{\Gamma}$ is such as in Part 3 of the theorem. 
This is obvious when $R$ is either trivial or 
a Galois field extension $L/K$ with $\Gamma=\mr{Gal}(L/K)$; notice from the preceding case that $\g_0(R)^{\Gamma}=\g_0(K)$
if $R$ is trivial. 
We may thus exclude these two cases. 

To describe $R$, note that $R$ is artinian as a ring, and $\Gamma$-simple in the sense
that it does not include any non-trivial $\Gamma$-stable ideal. Since the action by $\Gamma$ on $R$
commutes with the $\d$-operator, $R$ is a module algebra over the $\c$-Hopf algebra 
$\c \Gamma \otimes_{\c} \c[\d]$, which is \emph{artinian simple} or \emph{AS} in the sense of \cite[Definition 11.6]{AMT};
see the original \cite[Definition 2.6]{AM} as an alternate. 
This $\c$-Hopf algebra is the group algebra $\c \Gamma$ tensored with the polynomial algebra $\c[\d]$ in which $\d$ is primitive. 
Choose arbitrary a maximal (or equally, minimal) ideal $\mathfrak{m}$ of $R$, and let $\Gamma'$ be the subgroup of $\Gamma$ 
consisting of all elements that stabilize $\mathfrak{m}$. 
By \cite[Proposition 11.5]{AMT}
we have 
\[
\text{(a)}\ \Gamma' \simeq \mathbb{Z}_2\quad \text{or}\quad \text{(b)}\ \Gamma' \simeq \mathbb{Z}_3,
\] 
with the extremal cases being excluded.
Moreover, there exists a $\dK$-field $L$ such that $R$ is naturally
isomorphic to the $(\d,\Gamma)$-$K$-algebra $\op{Map}(\Gamma'\backslash \Gamma, L)$ consisting of all maps from the
set of right cosets $\Gamma'\backslash \Gamma$ to $L$. This $\op{Map}(\Gamma'\backslash \Gamma, L)$ is naturally isomorphic
to the direct product of $[\Gamma :\Gamma']$-copies of $L$, as $\dK$-algebra, and possesses the $\Gamma$-action presented below. 
Suppose that $\mathbb{Z}_2$ is an arbitrarily chosen 
subgroup of $\Gamma$ of order 2, and
$\mathbb{Z}_3$ is the unique subgroup of $\Gamma$ of order 3,
so that we have $\Gamma=\mathbb{Z}_3\rtimes \mathbb{Z}_2$. 

Case (a).\ We may suppose $\Gamma'=\mathbb{Z}_2$ (see \cite[Proposition 11.5 (1)]{AMT})
and $\Gamma'\backslash \Gamma=\mathbb{Z}_3$. If $\gamma\in \Gamma$, $\gamma'\in \mathbb{Z}_3$ 
and $f \in \op{Map}(\mathbb{Z}_3, L)$, then we have
\[
\gamma f : \gamma' \mapsto
\begin{cases}
f(\gamma\gamma'), &\text{if}\ \gamma \in \mathbb{Z}_3\\
\gamma f(\gamma'^{-1}), &\text{if}\ 0\ne \gamma \in \mathbb{Z}_2.
\end{cases}
\]

Case (b).\ We have $\Gamma'=\mathbb{Z}_3$, and 
we may suppose $\Gamma'\backslash \Gamma=\mathbb{Z}_2$. If $\gamma\in \Gamma$,
$\gamma'\in \mathbb{Z}_2$ 
and $f \in \op{Map}(\mathbb{Z}_2, L)$, then we have
\[
\gamma f : \gamma' \mapsto
\begin{cases}
f(\gamma\gamma'), &\text{if}\ \gamma \in \mathbb{Z}_2;\\
\gamma f(\gamma'), &\text{if}\ \gamma \in \mathbb{Z}_3,\ \gamma'=0~\text{in}~\mathbb{Z}_2;\\
\gamma^{-1}f(\gamma'), &\text{if}\ \gamma \in \mathbb{Z}_3,\ \gamma'\ne 0~\text{in}~\mathbb{Z}_2.
\end{cases}
\]

In either case, since $R$ is right $(K\Gamma)^*$-Galois, $\Gamma'$ must act non-trivially on $L$, so that
$L/K$ is a Galois field extension with $\Gamma'=\op{Gal}(L/K)$. Conversely, if $L/K$ is such, 
then $R$ is seen to be a right $(K\Gamma)^*$-Galois $\dK$-algebra, being split by $L$. 
Moreover, $\g_0(R)$ is naturally isomorphic to the $(\d,\Gamma)\text{-}R$-Lie algebra $\op{Map}(\Gamma'\backslash\Gamma, \g_0(L))$
equipped with the obviously induced structure. We see
\begin{align*}
\g_0(R)^{\Gamma}&\simeq\op{Map}(\Gamma'\backslash\Gamma, \g_0(L))^{\Gamma}
=\big( \op{Map}( \Gamma'\backslash\Gamma, \g_0(L) )^{\mathbb{Z}_3} \big)^{\mathbb{Z}_2}\\
&=\begin{cases}
\{\text{all constant maps}~\mathbb{Z}_3\to \g_0(L)\}^{\mathbb{Z}_2} & \text{in Case (a)}\\
\, \op{Map}(\mathbb{Z}_2, \g_0(L)^{\mathbb{Z}_3})^{\mathbb{Z}_2} & \text{in Case (b)}
\end{cases}\\
&=\g_0(L)^{\Gamma'},
\end{align*}
which completes the proof.
\medskip


\subsection{Explicit non-trivial forms}\label{subsec:explicit_form}
Let us describe explicitly (up to quasi-isomorphism) the non-trivial twisted forms of $\g_0(K)$ listed
in (ii) of Part 2 and (ii)--(iv) of Part 3 of the theorem, separately for four types. For all those, quadratic 
field extensions are needed. Such an extension $L/K$ is of the form 
\[ L=K(\sqrt{\alpha})=\{ a+b\sqrt{\alpha}\mid a, b \in K\},  \]
where $\alpha \in K^{\times}\setminus (K^{\times})^2$. 
The generator of $\mr{Gal}(L/K)\, (=\mathbb{Z}_2)$ sends each element $x=a+b\sqrt{\alpha}$ to
\begin{equation}\label{eq:overline}
\overline{x}:= a-b\sqrt{\alpha}. 
\end{equation}
We will use this symbol $\overline{x}$, regardless of $\alpha$.

\subsubsection{{\bf Type} $\mathrm{A}_{\ell}\ (\ell \ge 2)$}\
We have $\g_0=\mathfrak{sl}_n(\c)$, where $n=\ell+1\ge 3$. 
The order $2$ outer-automorphism is conjugate to
$X \mapsto -{}^tX$. 
For a quadratic extension field $L=K(\sqrt{\alpha})$ over $K$ as above, the generator
of $\Gamma\, (=\mathbb{Z}_2)$ may supposed to act on $\g_0(L)$ by
$X=\big(x_{ij}\big)_{i,j} \mapsto -{}^t\overline{X}=\big(-\overline{x}_{ji}\big)_{i,j}$; see \cite[Chapter IX, Theorem 5]{J}. 
We see
\[
\g_0(L)^{\Gamma}=\mathfrak{o}_n(K)~\oplus \sqrt{\alpha}\big(\op{Sym}_n(K)\cap \mathfrak{sl}_n(K)\big),
\]
where $\op{Sym}_n(K)$ (resp., $\mathfrak{o}_n(K)$) denotes the $K$-subspace of $\mathfrak{gl}_n(L)$ consisting of all 
matrices $X$ with entries in $K$ that are symmetric (resp., skew-symmetric, ${}^tX=-X$).

\subsubsection{{\bf Type} $\mathrm{D}_{\ell}\ (\ell \ge 5)$}\
Let $m=2\ell$. We have
$
\g_0=\mathfrak{o}_{m}(\c),
$
which consists of all skew-symmetric $m\times m$ complex matrices. 
The order $2$ outer-automorphism is conjugate to $X\mapsto DXD$, where $D=\op{diag}(-1,1, \dots,1)$; 
see \cite[Chapter IX, Theorem 6]{J}. 
For a quadratic extension field $L=K(\sqrt{\alpha})$ over $K$ as above, we see
\begin{equation}\label{eq:DGmma_inv}
\g_0(L)^{\Gamma}=
\left\{ \, \begin{pmatrix} 0& -\sqrt{\alpha} \, {}^tX \\ \sqrt{\alpha} \, X  &Y \end{pmatrix}\, \biggm|
X\in K^{m-1},\ Y \in \mathfrak{o}_{m-1}(K) \, \right\},
\end{equation} 
where by writing $X \in K^{m-1}$, we mean that $X$ is an 
$(m-1)$-columned vector with entries in $K$. 

\subsubsection{{\bf Type} $\mathrm{E}_6$}\
Here we follow Jacobson \cite[Section~7]{Jacobson} for the construction.
Let $\mathfrak{J}$ be the exceptional central simple Jordan algebra over $\c$, and let $\mathfrak{J}^+$ denote
the subspace of $\mathfrak{J}$ which consists of the elements $a$ with trace zero, $T(a)=0$. 
We have the general linear complex Lie algebra $\mathfrak{gl}(\mathfrak{J})$ on the $\c$-vector space $\mathfrak{J}$.
Given an element $a \in \mathfrak{J}^+$, we have an element $R_a \in \mathfrak{gl}(\mathfrak{J})$ given by
$R_a(x)= xa\, (=ax)$, $x \in \mathfrak{J}$. Let $R_{\mathfrak{J}^+}$ be the subspace of $\mathfrak{gl}(\mathfrak{J})$
which consists of all $R_a$, $a\in \mathfrak{J}^+$. 
The complex simple Lie algebra $\g_0$ of type $\mathrm{E}_6$ is 
the Lie subalgebra of $\mathfrak{gl}(\mathfrak{J})$ generated by $R_{\mathfrak{J}^+}$. We have
\[ \g_0= R_{\mathfrak{J}^+} \oplus~\mathfrak{f}_0, \]
where we set $\mathfrak{f}_0:=[R_{\mathfrak{J}^+},R_{\mathfrak{J}^+}]$; this is a Lie subalgebra of $\g_0$ such that
$[R_{\mathfrak{J}^+},\mathfrak{f}_0]=R_{\mathfrak{J}^+}$,
and is in fact the complex simple Lie algebra
of type $\mathrm{E}_6$. The order $2$ outer-automorphism of $\g_0$ is conjugate to $X\mapsto -X^*$,
where $X^*$ denotes the operator adjoint to $X$ with respect to the trace form $(a, b)\mapsto T(ab)$. 
More explicitly this is given by
\[
X \mapsto 
\begin{cases} 
-X & \text{if}\ X\in R_{\mathfrak{J}^+};\\
\ X & \text{if}\ X \in \mathfrak{f}_0.
\end{cases}
\]
Therefore, we have
\[ \g_0(L)^{\Gamma}= (R_{\mathfrak{J}^+}\, \ot_{\c} K\sqrt{\alpha}) \oplus~\mathfrak{f}_0(K). \]

\subsubsection{{\bf Type} $\mathrm{D}_4$}\ 
The complex Lie algebra $\g_0$ is the Lie algebra
$\mathfrak{o}_{8}(\c)$ of skew-symmetric $8 \times 8$ complex matrices. 
We follow \'{E}.~Cartan \cite{Cartan} for the explicit description 
of outer-automorphisms.
We discuss for each group action,
separately as in Part 3 of the theorem. 

(ii)\ {\bf Action by} $\mathbb{Z}_2$. The argument above for $\mathrm{D}_{\ell}$ $(\ell\ge 5)$ works
for $\ell=4$, as well, so that $\g_0(L)^{\mathbb{Z}_2}$ is given by the right-hand side of 
\eqref{eq:DGmma_inv} with $m=8$.

(iii)\ {\bf Action by} $\mathbb{Z}_3$. Choose a generator $\sigma$ of the group. 
The relevant Galois extension is a cubic one, and it is of the form 
$L=K(\sqrt[ 3 \, ]{\beta})$, where $\beta \in K^{\times}\setminus (K^{\times})^3$. 
The generator $\sigma$ 
acts on $L$ so that $\sqrt[ 3 \, ]{\beta}\mapsto \omega \sqrt[ 3 \, ]{\beta}$, where $\omega$ is a primitive 3rd root of $1$.

We suppose that the rows and the columns of matrices in $\g_0\, (=\mathfrak{o}_8(\c))$ are indexed by the eight integers
$0,1,\cdots, 7$. Given a matrix $X=\big( x_{ij} \big)_{0\le i,j\le 7}$ in $\g_0$, we define seven vectors in $\c^4$ by
\begin{equation}\label{eq:seven_vectors}
X_i= {}^t\big(x_{0,i},\ x_{i+1,i+5},\ x_{i+4,i+6},\ x_{i+2,i+3}\big),\quad i=1,2,\dots,7,
\end{equation}
where the index $i+p$ greater than 7 is understood to be $i+p-7$; e.g., the third entry $x_{i+4,i+6}$ in $i=3$
in understood to be $x_{7,2}\, (=-x_{2,7})$. Then every matrix $X$ as above is uniquely determined by these seven vectors.
This holds when $\g_0$ is replaced by its base extension. We will use in \eqref{Eq1}--\eqref{Eq3}
the notation $X_i$ for the seven vectors which are associated 
as above with a matrix $X$ in such a base extension.

The action by $\mathbb{Z}_3$ on $\g_0$ is (up to conjugation)  
such that $\sigma$ acts on the vectors above
as the $\c$-linear automorphisms given by the matrix
\begin{equation}\label{eq:Ess}
S=\frac{1}{2}
\begin{pmatrix}
-1&-1&-1&-1\\ 
\ 1&\ 1&-1&-1\\
\ 1&-1&\ 1&-1\\
\ 1&-1&-1&\ 1
\end{pmatrix},
\end{equation}
which is seen to have $1,1,\omega$ and $\omega^2$ as eigen-values; see \cite[Section 4]{Cartan}. 
Set $\sqrt{-3}:=1+2\omega$, a square root of $-3$. Then we have the eigen-vectors
\begin{equation}\label{eq:eigen-vectors}
\boldsymbol{v}^{(1)}_1=\begin{pmatrix} 0\\ 1\\-1\\0 \end{pmatrix},\
\boldsymbol{v}^{(2)}_1=\begin{pmatrix} 0\\ 1 \\ 0 \\ -1\end{pmatrix},\
\boldsymbol{v}_{\omega}=\begin{pmatrix} \sqrt{-3} \\ 1\\ 1\\1 \end{pmatrix},\
\boldsymbol{v}_{\omega^2}=\begin{pmatrix} -\sqrt{-3}\\ 1\\ 1\\ 1 \end{pmatrix}
\end{equation}
of $S$ which are associated with $1$, $1$, $\omega$ and $\omega^2$, respectively; these form a basis of $\c^4$. 
Let $L^4$ denote the $L$-vector space of all $4$-columned vectors with entries in $L$.
Define a 4-dimensional $K$-subspace of $L^4$ by
\begin{equation}\label{eq:Xi}
\Xi_{L/K}=
K \boldsymbol{v}^{(1)}_1+K\boldsymbol{v}^{(2)}_1
+K\sqrt[ 3 \, ]{\beta}\, \boldsymbol{v}_{\omega^2}
+K(\sqrt[ 3 \, ]{\beta})^2 \boldsymbol{v}_{\omega}.
\end{equation}
We see now easily
\begin{equation}\label{Eq1}
\g_0(L)^{\Gamma} =
\{ \,  X \in \g_0(L)\, \mid \,
X_i
\in \Xi_{L/K},\ 1\le i \le 7 \,
\}. 
\end{equation}

(iv)\ {\bf Action by} $\Gamma\, (=\mathfrak{S}_3)$.
The relevant Galois extension is described by the following.

\begin{lemma}\label{lem:GaloisS_3}
A Galois extension field over $K$ with Galois group $\Gamma\, (=\mathfrak{S}_3)$ is the same as a field $L$ 
of the form $L=K(\sqrt{\alpha}, \sqrt[ 3 \, ]{\beta})$, where
\begin{itemize}
\item[\textup{(a)}] $\alpha \in K^{\times}\setminus (K^{\times})^2$, so that $K(\sqrt{\alpha})/K$ is a quadratic field extension, 
\item[\textup{(b)}] $\beta \in M^{\times}\setminus (M^{\times})^3$, where $M=K(\sqrt{\alpha})$, and
\item[\textup{(c)}] $\beta\overline{\beta}\in (K^{\times})^3$, where $\overline{\beta}$ is such as given by \eqref{eq:overline}.
\end{itemize}
For such an $L$, we have
\begin{itemize}
\item[\textup{(x)}]
an order $3$ element $\sigma$ and an order $2$ element $\tau$ of $\Gamma$, which necessarily generate $\Gamma$, 
satisfying $\sigma\tau=\tau\sigma^2$, 
\item[\textup{(y)}]
a primitive 3rd root $\omega$ of $1$, and 
\item[\textup{(z)}]
an element $\gamma$ of $K^{\times}$ such that $\gamma^3=\beta\overline{\beta}$ (see \textup{(c)} above), 
\end{itemize}
with which the action by $\Gamma$ is presented as
\begin{equation*}
\sigma: \sqrt{\alpha}\mapsto \sqrt{\alpha},\ \sqrt[ 3 \, ]{\beta}\mapsto \omega \sqrt[ 3 \, ]{\beta};\quad
\tau: 
\sqrt{\alpha}\mapsto -\sqrt{\alpha}\, (=\overline{\sqrt{\alpha}}),\ 
\sqrt[ 3 \, ]{\beta}\mapsto\frac{\gamma}{\sqrt[ 3 \, ]{\beta}}.
\end{equation*}
\end{lemma}
\pf
Given $\beta$ such as in (b), we have a cubic extension $M(\sqrt[ 3 \, ]{\beta})/M$. 
One sees that $\beta\overline{\beta}\in (M^{\times})^3$ if and only if $M(\sqrt[ 3 \, ]{\beta})=M(\sqrt[ 3 \, ]{\overline{\beta}})$. 
If $\gamma \in M^{\times}$ and $\gamma^3=\beta\overline{\beta}$, then $\gamma/\sqrt[ 3 \, ]{\beta}$ is a 3rd root 
of $\overline{\beta}$. A point is only to see that $\sqrt[ 3 \, ]{\beta}\mapsto\gamma/\sqrt[ 3 \, ]{\beta}$ gives
an involution which extends $M \to M, x \mapsto \overline{x}$ if and only if $\gamma \in K^{\times}$. 
\epf

\begin{example}\label{example:Galois_extension}
Recall $K=\C$. One can prove that
\[
\alpha =1-t^3\quad \text{and}\quad  \beta=1+\sqrt{1-t^3} 
\]
satisfy the conditions above.
Indeed, a point is to prove $\beta \notin (M^{\times})^3$, reducing it to show directly that there is no triple
of polynomials $a, b, c$ in $\c[t]$ with $c$ monic, such that
\[
a^3 + 3ab^2\alpha=c^3,\quad 3a^2b+b^3\alpha=c^3. 
\]
The result shows that there exists a Galois extension $L/K$ with $\op{Gal}(L/K)=\mathfrak{S}_3$.
\end{example}

Let 
$L=K(\sqrt{\alpha},\ \sqrt[ 3 \, ]{\beta}),\ M=K(\sqrt{\alpha})$,\ 
$\sigma$, $\tau$, $\omega$ and $\gamma$ be as in Lemma \ref{lem:GaloisS_3}. 
Recall from the proof of the lemma that $\gamma/\sqrt[ 3 \, ]{\beta}$ is a 3rd root 
of $\overline{\beta}$, and denote it by $\sqrt[ 3 \, ]{\overline{\beta}}$, so that one has
\[
\tau(\sqrt[ 3 \, ]{\beta})=\sqrt[ 3 \, ]{\overline{\beta}},\quad
\tau(\sqrt[ 3 \, ]{\overline{\beta}})=\sqrt[ 3 \, ]{\beta} .
\]

The action by $\Gamma$ on $\g_0\, (=\mathfrak{o}_8(\c))$ is (up to conjugation)  
such that the generators $\sigma$ and $\tau$ act on the seven vectors in \eqref{eq:seven_vectors} 
as the $K$-linear automorphisms given by the matrix $S$ in \eqref{eq:Ess} and the diagonal
matrix $D=\op{diag}(-1,1,1,1)$, respectively. The latter action by $\tau$ on $\g_0$ 
coincides with the above mentioned outer-automorphism $X\mapsto DXD$ for type $\mathrm{D}_{\ell}$,
when $\ell=4$. 

Note $L=M(\sqrt[ 3 \, ]{\beta})$, and
apply the previous result for the action by $\mathbb{Z}_3$ to the action by $\langle \sigma \rangle$ (on $L/M$). 
Then, by using the $M$-subspace $\Xi_{L/M}$ of $L^4$ defined by \eqref{eq:Xi} (modified into the present situation), we have
\begin{equation}\label{Eq2}
\g_0(L)^{\langle \sigma \rangle} = 
\{ \, X \in \g_0(L)\, \mid \,
X_i
\in \Xi_{L/M},\ 1\le i \le 7 \,
\}. 
\end{equation}
By using the vectors given in \eqref{eq:eigen-vectors}
we define a 4-dimensional $K$-subspace of $L^4$ by 
\begin{align*}
\Theta_{L/K}=
K\boldsymbol{v}^{(1)}_1+K\boldsymbol{v}^{(2)}_1
+K\Big(\sqrt[ 3 \, ]{\beta}\, \boldsymbol{v}_{\omega^2}+\sqrt[ 3 \, ]{\overline{\beta}}\, \boldsymbol{v}_{\omega}\Big)~\\ +
K\sqrt{\alpha}\Big(\sqrt[ 3 \, ]{\beta}\, \boldsymbol{v}_{\omega^2}-\sqrt[ 3 \, ]{\overline{\beta}}\, \boldsymbol{v}_{\omega}\Big).
\end{align*}
We see now easily
\begin{equation}\label{Eq3}
\g_0(L)^{\Gamma} = (\g_0(L)^{\langle \sigma \rangle})^{\langle \tau \rangle} =
\{ \, X \in \g_0(L)\, \mid \,
X_i
\in \Theta_{L/K},\ 1\le i \le 7
\, \}. 
\end{equation}

\section*{Acknowledgments}
The authors thank Professor Arturo Pianzola who kindly brought to them an interesting problem, informing of Steinberg's Theorem,
and gave helpful comments during their revision of the manuscript. 
They also thank the referees for helpful suggestions and valuable comments.


\begin{thebibliography}{99}

\bibitem{AM} 
K.~Amano,\ A.~Masuoka,\ 
\emph{Picard-Vessiot extensions of artinian simple module algebras},\
J.\ Algebra {\bf 285} (2005),\ 743--767. 


\bibitem{AMT} 
K.~Amano,\ A.~Masuoka,\ M.~Takeuchi,
\emph{Hopf algebraic approach to Picard-Vessiot theory},\
in: M.~Hazewinkel (ed.), \emph{Handbook of Algebra}, 
Vol. 6, Elsevier, 
North-Holland, 2009, pp.~127--171. 

\bibitem{B1} 
J.~Bausch,\ 
\emph{Automorphismes des alg\`{e}bres de Kac-Moody affines} (French)\
[Automorphisms of affine Kac-Moody algebras], C. R. Acad. Sci. Paris S\`{e}r. I Math.
{\bf 302} (1986), no. 11, 409--412.

\bibitem{B2}
J.~Bausch,\ 
\emph{\'{E}tude et classification des automorphismes d'ordre fini et de premi\`{e}re
esp\`{e}ce des alg\`{e}bres de Kac-Moody affines} (French)\ [Study and classification of
automorphisms of finite order and the first kind of affine Kac-Moody algebras] with
appendices by G.~Rousseau,\ 
Inst. \`{E}lie Cartan 11, Alg\`{e}bres de Kac-Moody affines, pp. 5--124, 
Univ. Nancy, Nancy, 1989.

\bibitem{BaR} 
J.~Bausch,\ G.~Rousseau,\
\emph{Involutions de premi\`{e}re esp\`{e}ce des alg\`{e}bres affines}
(French), [Involutions of the first kind of affine algebras],  
Inst. \`{E}lie Cartan 11, Alg\`{e}bres de Kac-Moody affines, pp. 125--139, 
Univ. Nancy, Nancy, 1989.

\bibitem{BR1} H.~Ben Messaoud,\ G.~Rousseau,\
\emph{Classification des formes r\`{e}elles presque
compactes des alg\`{e}bres de Kac-Moody affines} (French summary),
[Classification of almost compact real forms of affine Kac-Moody algebras],\ 
J.\ Algebra {\bf 267} (2003), no. 2, 443--513.

\bibitem{BR2} H.~Ben Messaoud,\ G.~Rousseau,\
\emph{Erratum to ``Classification of almost compact
real forms of affine Kac-Moody algebras"},\ 
J.\ Algebra {\bf 279} (2004), no. 2, 850--851.


\bibitem{Cartan}
\'{E}.~Cartan,\
\emph{Le principe de dualit\'{e} et la the\'{o}rie des groupes simples et semi-simples},\
Bull. Sci. Math. {\bf 49} (1925), 361--374; \emph{Oeuvres compl\`{e}tes}, Partie I, 
\emph{Groupes de Lie},
\'{E}ditions du CNRS, 1984, pp. 555--568.



\bibitem{J} N.~Jacobson,\
\emph{Lie Algebras},\ 
Interscience Publ., New York, 1962. 

\bibitem{Jacobson} N.~Jacobson,\
\emph{Exceptional Lie Algebras},\ 
Lec. Notes in Pure and Appl. Math., Vol.~1, Marcel Dekker, Inc., New York, 1971. 

\bibitem{Katz} N.~M.~Katz,\ 
\emph{A conjecture in the arithmetic theory of differential equations},\
Bull. Soc. Math. France {\bf 110} (1982), no. 2, 203--239; \emph{Corrections},\ 
Bull. Soc. Math. France {\bf 110} (1982), no. 3,
347--348.


\bibitem{KMRT} M-A.~Knus,\ A.~Merkurjev,\ M.~Rost,\ J-P.~Tignol,\ 
\emph{The Book of Involutions},\
Colloquium Publications Vol.~44, Amer. Math. Soc., Providence, 1998. 


\bibitem{K1} E.~R.~Kolchin,\ 
\emph{Differential Algebra and Algebraic Groups},\
Pure an Applied Mathematics, Vol.~54, Academic Press, New York, 1973. 


\bibitem{Mo} S.~Montgomery,\
\emph{Hopf Algebras and Their Actions on Rings},\
CBMS Conf.\ Series in Math., Vol.~82, Amer.\ Math.\ Soc., Providence, RI, 1993.

\bibitem{OV} A.~L.~Onishchik,\ E.~B.~Vinberg,\
\emph{Lie Groups and Algebraic Groups},\ 
Sprinfer-Verlag, Heidelberg, 1990.

\bibitem{P} A.~Pianzola,\
\emph{Vanishing of $H^1$ for Dedekind rings and applications to loop algebras},\ 
C. R. Acad. Sci. Paris, Ser. I, {\bf 340} (2005), 633--638.

\bibitem{VS} 
M.~van~der Put,\ M.~F.~Singer,\
\emph{Galois Theory of Linear Differential Equations},\
Grundlehren der mathematischen Wissenschaften, Vol.~328, Springer-Verlag, Berlin/Heidelberg, 2003. 

\bibitem{Sc} P.~Schauenburg,\
\emph{Hopf bigalois extensions},\ 
Comm. Algebra {\bf 24} (1996), no.~12, 3797--3825.


\bibitem{Ser1}
V.~Serganova,\
\emph{Classification of simple real Lie superalgebras and symmetric
superspaces}, Funktsional. Anal. i Prilozhen {\bf 17} (1983), no. 3, 46--54 (Russian);
English translation: Functional Anal. Appl. {\bf 17} (1983), no. 3, 200--207.

\bibitem{Ser2}
V.~Serganova,\
\emph{Automorphisms of Lie superalgebras of string theories},
Funktsional. Anal. i Prilozhen {\bf 19} (1985), no. 3, 75--76 (Russian); 
English translation: Functional Anal. Appl. {\bf 19} (1985), no. 3, 226--228.


\bibitem{Se} J.-P.~Serre,\
\emph{Galois Cohomology},\
Springer-Verlag, Berlin/Heidelberg, 1997.


\bibitem{W}
W.~C.~Waterhause,\
\emph{Introduction to Affine Group Schemes},
Graduate Texts in Mathematics, Vol.~66, 
Springer-Verlag, New York, 1979. 

\end{thebibliography}
\end{document}